\documentclass{article}

\usepackage[english]{babel}
\usepackage{amsmath}
\usepackage{amssymb}
\usepackage{marvosym}
\usepackage{amsfonts}
\usepackage{array}
\usepackage[dvips]{graphicx}
\include{bm}
\usepackage{epsfig}
\usepackage{float}
\usepackage{comment}

\usepackage{amsmath,amssymb,amsfonts,enumerate,rotate,multirow}
\usepackage{bm,latexsym,mathrsfs}
\usepackage{epsfig,graphicx}
\usepackage[usenames]{color}
\usepackage{theorem}
\newtheorem{thm}{Theorem}[section]

\newtheorem{lem}[thm]{Lemma}

\theoremstyle{remark}
\newtheorem{rem}[thm]{Remark}

\numberwithin{equation}{section}

 \newcommand{\ind}{\mathbf{1}}

 \renewcommand{\tilde}{\widetilde}

 \def\proof{\noindent {\bf Proof.}\ }
\def\endproof{{\mbox{}\nolinebreak\hfill\rule{2mm}{2mm}\par\medbreak} }

\begin{document}

\title{\textbf{Dynamics of an Adaptive\\ Randomly Reinforced Urn}}
\author{Giacomo Aletti$^a$,
Andrea Ghiglietti$^a$ and
Anand N. Vidyashankar$^b$}

\maketitle
\begin{center}
\small $^a$ADAMSS center \& Department of Mathematics, Universit\`{a} degli Studi di Milano, Milan, Italy\\
\small $^b$Department of Statistics, George Mason University, Fairfax (VA), USA
\end{center}

\abstract{Adaptive randomly reinforced urn (ARRU) is a two-color urn
model where the updating process is defined by a sequence of non-negative
random vectors $\{(D_{1,n}, D_{2,n});n\geq1\}$ and randomly evolving thresholds
which utilize accruing statistical information for the updates. Let
$m_1=E[D_{1,n}]$ and $m_2=E[D_{2,n}]$. Motivated by applications, in this paper
we undertake a detailed study of the dynamics of the ARRU model. First,
for the case $m_1 \neq m_2$, we establish $L_1$ bounds on the increments of
the urn proportion at fixed and increasing times under very weak assumptions on the
random threshold sequence. As a consequence, we deduce weak
consistency of the evolving urn proportions. Second, under
slightly stronger conditions, we establish the strong consistency of the
urn proportions for all finite values of $m_1$ and $m_2$. Specifically,
we show that when $m_1=m_2$ the proportion converges to a non-degenerate
random variable. Third, we establish the asymptotic distribution, after
appropriate centering and scaling, of the proportion of sampled balls in the
case $m_1=m_2$. In the process, we settle the issue of asymptotic
distribution of the number of sampled balls for a randomly reinforced urn (RRU).
To address the technical issues, we establish results on the harmonic
moments of the total number of balls in the urn at different times under
very weak conditions, which is of independent interest.}

\noindent
{\bf Keywords}:
generalized P\'olya urn, reinforced processes, strong and weak consistency, central limit theorems, crossing times, harmonic moments.

\noindent
{\bf MSC Subject Classification}: 60F05, 60F15, 97K50.

\section{Introduction}   \label{section_introduction}

In recent years, randomly reinforced urn (RRU) has been investigated in statistical and probability literature as a model for clinical trial design,
computer experiments and in the context of vertex reinforced random walk (see~\cite{Hu.et.al.06,Mahmoud.08,Pemantle.et.al.99}).
Introduction of accruing information in the implementation of these urn models in practice, leads to an adaptive randomly reinforced urn (ARRU).
In this paper, we study the properties concerning the urn composition of an ARRU.
We now turn to a precise description of the ARRU.

A randomly reinforced urn (RRU) model (see~\cite{Muliere.et.al.06}) is characterized by a pair $(Y_{1,n},Y_{2,n})$ of real random variables representing
the number of balls of two colors, red and white.
The process is described as follows: at time $n=0$, the process starts with $(y_{1,0},y_{2,0})$ balls.
A ball is drawn at random. If the color is red, the ball is returned to the urn along with the random numbers $D_{1,1}$ of red balls;
otherwise, the ball is returned to the urn along with the random numbers $D_{2,1}$ of white balls.
Let $Y_{1,1}=y_{1,0}+D_{1,1}$ and $Y_{2,1}=y_{2,0}$ denote the urn composition when the sampled ball is red;
similarly, let $Y_{1,1}=y_{1,0}$ and $Y_{2,1}=y_{2,0}+D_{2,1}$ denote the urn composition when the sampled ball is white.
The process is repeated yielding the collection $\{(Y_{1,n},Y_{2,n});n\geq1\}$.
The quantities $\{D_{1,n};n\geq1\}$ and $\{D_{2,n};n\geq1\}$ are independent collections of independent and identically distributed (i.i.d.)
non-negative random variables.
Hence, an RRU model is characterized by the \emph{replacement matrix}
\[\mathbf{D}_{n}=\left[
\begin{array}{ll}
D_{1,n}\ &\ 0\\
0\ &\ D_{2,n}\\
\end{array}
\right].\]
Let $m_1:=\bm{E}[D_{1,n}]$ and $m_{2}:=\bm{E}[D_{2,n}]$.
The asymptotic properties of the urn composition in the above model were investigated by Muliere \textit{et al.} (see~\cite{Muliere.et.al.06})
and Aletti \textit{et al.} (see~\cite{Aletti.et.al.09});
specifically, they established that
\begin{equation}\label{eq:RRU}
Z_n\ =\ \frac{Y_{1,n}}{Y_{1,n}+Y_{2,n}}\ \stackrel{a.s.}{\rightarrow}\
\begin{cases}
1\ &\text{if } m_1> m_2,\\
Z_{\infty}\ &\text{if } m_1= m_2,\\
0\ &\text{if } m_1< m_2,
\end{cases}\end{equation}
where $\stackrel{a.s.}{\rightarrow}$ stands for almost sure convergence and $Z_{\infty}$ is a non-degenerate random variable supported on $(0,1)$.
The properties of the distribution of $Z_{\infty}$ were studied in Aletti \textit{et al.} (see~\cite{Aletti.et.al.09,Aletti.et.al.12}).
Specifically, it is shown in Aletti \textit{et al.} (see~\cite{Aletti.et.al.09}) that when $m_1=m_2$,
$\bm{P}(Z_{\infty}=x)=0$ for any $x\in[0,1]$.
Denoting $\{(N_{1,n},N_{2,n});n\geq1\}$ the number of balls of red and white colors sampled from the urn,
one can deduce from  (1.1) that $N_{1,n}/n$ converges to the same limit as $Z_n$.

Notice that the limit of the RRU in~\eqref{eq:RRU} is always 1 or 0 when $m_1\neq m_2$,
and the rate of convergence and the limit distribution has been established in May and Flournoy (2009) (see~\cite{Flournoy.et.al.09}).
However, motivated by applications in clinical trials (see~\cite{Hu.et.al.06}),
it is common to target a specific value $\rho\in(0,1)$.
This was achieved in Aletti \textit{et al.} (see~\cite{Aletti.et.al.13}),
where the modified randomly reinforced urn (MRRU) model was introduced.
The MRRU model is an RRU model with two fixed thresholds $0<\rho_2\leq \rho_1<1$,
such that if $Z_n<\rho_2$, no white balls are replaced in urn,
while if $Z_n>\rho_1$, no red balls are replaced in the urn.
The replacement matrix in this case is
\[\mathbf{D}_{n}=\left[
\begin{array}{ll}
D_{1,n}\cdot\ind_{\{Z_{n-1}\leq\rho_1\}}\ &\ 0 \\
0\ &\ D_{2,n}\cdot\ind_{\{Z_{n-1}\geq\rho_2\}}\\
\end{array}
\right].
\]
The strong consistency in the case $m_1\neq m_2$ was established in Aletti \textit{et al.} (see~\cite{Aletti.et.al.13}); i.e. they showed that
\begin{equation}\label{eq:MRRU}
Z_n\ \stackrel{a.s.}{\rightarrow}\
\begin{cases}
\rho_1\ &\text{if } m_1> m_2,\\
\rho_2\ &\text{if } m_1< m_2.
\end{cases}
\end{equation}
A second order result for $Z_n$, namely the asymptotic distribution of $Z_n$
after appropriate centering, was derived in Ghiglietti \textit{et al.} (see~\cite{Ghiglietti.et.al.12}).

In applications, especially in clinical trials (see~\cite{Hu.et.al.06}),
$\rho_1$ and $\rho_2$ are unknown and depend on the parameters of the distributions of $D_{1,1}$ and $D_{2,1}$.
Let $\mathcal{F}_{n-1}$ be the $\sigma$-algebra generated by the information up to time $n-1$ and
let $\hat{\rho}_{1,n-1}$ and $\hat{\rho}_{2,n-1}$ be two random variables that are $\mathcal{F}_{n-1}$-measurable.
Ghiglietti \textit{et al.} proposed in~\cite{Ghiglietti.et.al.15}
an adaptive randomly reinforced urn model that uses accruing information to construct random thresholds $\hat{\rho}_{1,n-1}$ and $\hat{\rho}_{2,n-1}$
which converge a.s. to specified targets $\rho_1$ and $\rho_2$. Thus, using the replacement matrix
\begin{equation}\label{eq:replacement_matrix}
\mathbf{D}_{n}=\left[
\begin{array}{ll}
D_{1,n}\cdot\ind_{\{Z_{n-1}\leq\hat{\rho}_{1,n-1}\}}\ &\ 0 \\
0\ &\ D_{2,n}\cdot\ind_{\{Z_{n-1}\geq\hat{\rho}_{2,n-1}\}}\\
\end{array}
\right],\end{equation}
an MRRU becomes an Adaptive Randomly Reinforced Urn (ARRU).
It is worth mentioning here that the random thresholds $\hat{\rho}_{1,n-1}$ and $\hat{\rho}_{2,n-1}$
depend on the \emph{adaptive} estimators of the parameters of the distributions of $D_{1,1}$ and $D_{2,1}$.

In a recent work,  Ghiglietti \textit{et al.} (see~\cite{Ghiglietti.et.al.15}) studied the asymptotic properties of an ARRU when $m_1 \ne m_2$ under strong conditions on the rate of convergence of the adaptive thresholds. Specifically, they established a strong consistency and asymptotic normality for the number of sampled balls under an exponential rate of convergence assumption on the adaptive thresholds. In this paper, first we establish that under very weak conditions, weak consistency of the proportion $Z_n$. This is achieved by providing useful and non-trivial $L_1$ bounds on (i) the increments of the distace $\Delta_n= |Z_n-\rho_1|$ (Theorem~\ref{thm:increments_Z}) and (ii) the increments of $\Delta_n$ at linearly increasing times (Theorem~\ref{thm:exponential_bounds_similar_MRRU_giu} and Theorem~\ref{thm:exponential_bounds_similar_MRRU_su}). These results provide insight into the dynamics of the ARRU  and are of independent interest. The proofs of these results need estimates on the harmonic moments of the total number of balls in the urn under weak assumptions on the thresholds. This result, of independent interest, is established in Theorem~\ref{thm:Y_bounded_and_diverge}. Second, we undertake a detailed analysis of the ARRU model when $m_1=m_2$. Specifically, we establish strong consistency of the proportion $Z_n$ and the limit distribution of the proportion of sampled balls for the ARRU. In the process, we also address the issue of limit distribution of the number of sampled balls from a randomly reinforced urn (RRU) thus settling one of the long-standing open problems in the field.

The rest of the paper is structured as follows: Section~\ref{section_the_model} contains the model, assumptions and main results;
Section~\ref{section_preliminary_results} is concerned with preliminary estimates and results on the urn process.
Sections~\ref{section_WLLN} and~\ref{section_SLLN} are concerned with the proofs of the consistency of the urn proportion
and Section~\ref{section_CLT} is concerned with the proof of the limit distribution of the proportion of sampled balls.

\section{Model assumptions, notation and main results}   \label{section_the_model}

We begin by describing our model precisely.
Let $\bm{\xi_1}=\left\{\xi_{1,n};n\geq1\right\}$ and $\bm{\xi_2}=\left\{\xi_{2,n};n\geq1\right\}$ be two sequences of i.i.d. random variables,
with probability distributions $\mu_1$ and $\mu_2$ respectively.
Without loss of generality (wlog), assume that the support of $\xi_{1,n}$ and $\xi_{2,n}$ is the same.
We denote it by $S$.
Consider an urn containing $y_{1,0}>0$ red balls and $y_{2,0}>0$ white balls, and define $y_0=y_{1,0}+y_{2,0}$ and $z_0 = y_0^{-1}y_{1,0}$.
In general, $y_{1,0}$ and $y_{2,0}$ may not assume integer values.
At time $n=1$, a ball is drawn at random from the urn and its color is observed.
Let the random variable $X_1$ be such that
\begin{equation*}
X_1\ =\
\begin{cases}
1\ &\text{if the extracted ball is red},\\
0\ &\text{if the extracted ball is white}.
\end{cases}\end{equation*}
We assume $X_1$ to be independent of the sequences $\bm{\xi_{1}}$ and $\bm{\xi_{2}}$.
To make this assumption more explicit, we define $X_1=\ind_{\{U_1\leq z_0\}}$,
where $U_1$ is a uniform random variable in (0,1) independent of $\bm{\xi_{1}}$ and $\bm{\xi_{2}}$.
Note that $X_1$ Bernoulli random variable with parameter $z_0$.

Let $\hat{\rho}_{1,0}$ and $\hat{\rho}_{2,0}$ be two random variables such that $\hat{\rho}_{1,0},\hat{\rho}_{2,0}\in[0,1]$ and $\hat{\rho}_{1,0}\geq\hat{\rho}_{2,0}$.
Let $u:S\rightarrow \left[a, b\right]$, $0< a \leq b < \infty$.
If $X_{1}=1$ and $z_0\leq\hat{\rho}_{1,0}$, we return the extracted ball to the urn together with $D_{1,1}=u\left(\xi_{1,1}\right)$ new red balls.
While, if $X_{1}=0$ and $z_0\geq\hat{\rho}_{2,0}$, we return it to the urn together with $D_{2,1}=u\left(\xi_{2,1}\right)$ new white balls.
If $X_{1}=1$ and $z_0>\hat{\rho}_{1,0}$, or if $X_{1}=0$ and $z_0<\hat{\rho}_{2,0}$, the urn composition is not modified.
To ease notation, let denote $w_{1,0}=\ind_{\{z_0\leq\hat{\rho}_{1,0}\}}$ and $w_{2,0}=\ind_{\{z_0\geq\hat{\rho}_{2,0}\}}$.
Formally, the extracted ball is always replaced in the urn together with
$$X_{1}D_{1,1}w_{1,0}+\left(1-X_{1}\right)D_{2,1}w_{2,0}$$
new balls of the same color;
now, the urn composition becomes
\[\left\{
\begin{array}{l}
Y_{1,1}=y_{1,0}+X_1D_{1,1}w_{1,0}\\
\\
Y_{2,1}=y_{2,0}+(1-X_1)D_{2,1}w_{2,0}.
\end{array}
\right.\]
Set $Y_{1}=Y_{1,1}+Y_{2,1}$ and $Z_{1}=Y^{-1}_1Y_{1,1}$.
Now, by iterating the above procedure we define
$\hat{\rho}_{1,1}$ and $\hat{\rho}_{2,1}$ to be two random variables, with $\hat{\rho}_{1,1},\hat{\rho}_{2,1}\in[0,1]$ and $\hat{\rho}_{1,1}\geq\hat{\rho}_{2,1}$ a.s., measurable with respect to the $\sigma$-algebra $\mathcal{F}_1=\sigma\left(\mathcal{G}_1,\varphi_1\right)$,
where $\mathcal{G}_1=\sigma\left(X_1,X_1\xi_{1,1}+(1-X_1)\xi_{2,1}\right)$ and $\varphi_1$ is a r.v. independent of $\mathcal{G}_1$.
Let $m_1 =\int u\left(y\right)\mu_1\left(dy\right)$ and $m_2 =\int u\left(y\right)\mu_2\left(dy\right)$ be the means of
$\{D_{1,n};n\geq1\}$ and $\{D_{2,n};n\geq1\}$ respectively.

The urn process is then repeated for all $n\geq1$.
Let $\hat{\rho}_{1,n}$ and $\hat{\rho}_{2,n}$ be two random variables with $\hat{\rho}_{1,n},\hat{\rho}_{2,n}\in\left(0,1\right)$ and $\hat{\rho}_{1,n}\geq\hat{\rho}_{2,n}$ a.s., measurable with respect to the $\sigma$-algebra $\mathcal{F}_n=\sigma\left(\mathcal{G}_n,\varphi_1,..,\varphi_n\right)$, where $$\mathcal{G}_n=\sigma\left(X_1,X_1\xi_{1,1}+\left(1-X_1\right)\xi_{2,1},...,X_n,X_n\xi_{1,n}+\left(1-X_n\right)\xi_{2,n}\right),$$
and $\varphi_n$ are a collection of r.v. independent of $\mathcal{G}_n$.
We will refer to $\hat{\rho}_{j,n}$ $j=1,2$ as threshold parameters.

At time $n+1$, a ball is extracted and let $X_{n+1}=1$ if the ball is red and $X_{n+1}=0$ otherwise.
Equivalently, we can define $X_{n+1}=\ind_{\{U_{n+1}\leq Z_n\}}$,
where $U_{n+1}$ is a uniform random variable in (0,1) independent of $\mathcal{F}_{n}$, $\bm{\xi_{1}}$ and $\bm{\xi_{2}}$.
Then, the ball is returned to the urn together with
$$X_{n+1}D_{1,n+1}W_{1,n}+\left(1-X_{n+1}\right)D_{2,n+1}W_{2,n}$$
balls of the same color, where $D_{1,n+1}=u\left(\xi_{1,n+1}\right)$, $D_{2,n+1}=u\left(\xi_{2,n+1}\right)$, $W_{1,n}=\ind_{\{Z_n\leq\hat{\rho}_{1,n}\}}$,
$W_{2,n}=\ind_{\{Z_n\geq\hat{\rho}_{2,n}\}}$ and $Z_{n+1}=Y_{1,{n+1}}/Y_{n+1}$ for any $n\geq1$,
where
\[\left\{
\begin{array}{l}
Y_{1,n+1}=y_{1,0}+\sum_{i=1}^{n+1}X_iD_{1,i}W_{1,i-1}\\
\\
Y_{2,n+1}=y_{2,0}+\sum_{i=1}^{n+1}\left(1-X_i\right)D_{2,i}W_{2,i-1}
\end{array}
\right.\]
and $Y_{n+1}=Y_{1,n+1}+Y_{2,n+1}$.
If $X_{n+1}=1$ and $Z_n>\hat{\rho}_{1,n}$, i.e. $W_{1,n}=0$, or if $X_{n+1}=0$ and $Z_n<\hat{\rho}_{2,n}$, i.e. $W_{2,n}=0$,
the urn composition does not change at time $n+1$.
Note that condition $\hat{\rho}_{1,n}\geq\hat{\rho}_{2,n}$ a.s., which implies $W_{1,n}+W_{2,n}\geq1$, ensures that the urn composition can change with positive probability for any $n\geq1$, since the replacement matrix is never a zero matrix.
Since, conditionally to the $\sigma$-algebra $\mathcal{F}_n$, $X_{n+1}$ is assumed to be independent of $\xi_{1},\xi_{2}$,
$X_{n+1}$ is Bernoulli distributed with parameter $Z_n$.

\subsection{Weak consistency of the urn composition}   \label{subsection_WLLN}

A particulary relevant result of this paper is concerned with the consistency of the urn proportion $Z_n$
when the random thresholds $\hat{\rho}_{1,n}$ and $\hat{\rho}_{2,n}$ converge in probability to some constants in $\rho_1,\rho_2\in(0,1)$.
To obtain this result, we need to assume that the thresholds sequence are bounded away from 0 and 1 with high probability,
which is expressed in the following condition:
there exist two constants $0<\rho_{\min}\leq\rho_{\max}<1$ and $0<c_{\rho}<\infty$ such that
\begin{equation}\label{ass:range_bounded}
\bm{P}\left(\rho_{\min}\leq\hat{\rho}_{2,n}\leq\hat{\rho}_{1,n}\leq\rho_{\max}\right)\ \geq\ 1-\exp\left(-c_{\rho}n\right)\,
\end{equation}
for large $n$.
Hence, we can establish the consistency result as follows
\begin{thm}\label{thm:convergence_probability}
Assume~\eqref{ass:range_bounded} and there exist two constant $\rho_1,\rho_2\in(0,1)$, with $\rho_1\geq\rho_2$, such that
\begin{equation}\label{eq:rho_convergence_probability}
\hat{\rho}_{1,n}\stackrel{p}{\rightarrow}\rho_1\ \qquad\qquad\qquad\ \hat{\rho}_{2,n}\stackrel{p}{\rightarrow}\rho_2.
\end{equation}
Then, when $m_1\neq m_2$,
\begin{equation}\label{eq:Z_convergence_probability}
Z_n\ \stackrel{p}{\rightarrow}\
\begin{cases}
\rho_1\ &\text{if } m_1> m_2,\\
\rho_2\ &\text{if } m_1< m_2.
\end{cases}\end{equation}
\end{thm}
We present the proof of Theorem~\ref{thm:convergence_probability} in Section~\ref{section_WLLN}.

\begin{rem}
The strong consistency of the urn proportion presented in Ghiglietti \textit{et al.} (see~\cite{Ghiglietti.et.al.15}), i.e. $\hat{\rho}_{1,n}\stackrel{a.s.}{\rightarrow}\rho_1$ implies
$Z_n\stackrel{a.s.}{\rightarrow}\rho_1$, may suggest to prove Theorem~\ref{thm:convergence_probability} by applying subsequence arguments.
Specifically, $Z_n\stackrel{p}{\rightarrow}\rho_1$ in~\eqref{eq:Z_convergence_probability}
implies that for any subsequence $\{n_k;k\geq1\}$ there exists
a further subsequence $\{n_{k_j};j\geq1\}$ such that $Z_{n_{k_j}}\stackrel{a.s.}{\rightarrow}\rho_1$.
Moreover, assumption $\hat{\rho}_{1,n}\stackrel{p}{\rightarrow}\rho_1$ in~\eqref{eq:rho_convergence_probability} guarantees the existence of
$\{n_{k_j};j\geq1\}$ such that $\hat{\rho}_{1,n_{k_j}}\stackrel{a.s.}{\rightarrow}\rho_1$.
Nevertheless, the strong consistency result in~\cite{Ghiglietti.et.al.15} does not prove that $Z_{n_{k_j}}\stackrel{a.s.}{\rightarrow}\rho_1$ with the only assumption that $\hat{\rho}_{1,n_{k_j}}\stackrel{a.s.}{\rightarrow}\rho_1$, because
this condition does not provide any information on the behavior of $\hat{\rho}_{1,i}$ at times $i\notin\{n_{k_j};j\geq1\}$.
Hence, the convergence of $\hat{\rho}_{1,n_{k_j}}$ would imply the convergence of $Z_{n_{k_j}}$
only if the urn composition was updated exclusively at times $\{n_{k_j};j\geq1\}$.
\end{rem}

\subsection{Strong consistency of the urn composition}   \label{subsection_SLLN}

The following theorem states the consistency of the urn proportion $Z_n$ for any values of $m_1$ and $m_2$,
when the random thresholds $\hat{\rho}_{1,n}$ and $\hat{\rho}_{2,n}$ converge with probability one.
\begin{thm}\label{thm:convergence_as}
Assume there exist two constant $\rho_1,\rho_2\in[0,1]$, with $\rho_1\geq\rho_2$, such that
\begin{equation}\label{eq:rho_convergence_as}
\hat{\rho}_{1,n}\stackrel{a.s.}{\rightarrow}\rho_1\ \qquad\qquad\qquad\ \hat{\rho}_{2,n}\stackrel{a.s.}{\rightarrow}\rho_2.
\end{equation}
Then,
\begin{equation}\label{eq:Z_convergence_as}
Z_n\ \stackrel{a.s.}{\rightarrow}\
\begin{cases}
\rho_1\ &\text{if } m_1> m_2,\\
Z_{\infty}\ &\text{if } m_1= m_2,\\
\rho_2\ &\text{if } m_1< m_2,
\end{cases}\end{equation}
where $Z_{\infty}$ is a random variable such that $\bm{P}(Z_{\infty}\in[\rho_2,\rho_1])=1$.
\end{thm}
We present the proof of Theorem~\ref{thm:convergence_as} in Section~\ref{section_SLLN}.
When the limit of the urn proportion is different from 1 or 0,
the following convergence result on the total number of balls to the smaller mean holds.
\begin{lem}\label{lem:harmonic_as_convergence}
Assume~\eqref{eq:rho_convergence_as} with $\rho_1>\rho_2$ and let $m^{*}=\min\{m_1,m_2\}$.
Then, on the set $\{\lim_{n\rightarrow\infty}Z_n\neq\{0,1\}\}$,
$$\frac{Y_n}{n}\ \stackrel{a.s.}{\rightarrow}\ m^{*}.$$
\end{lem}
The above lemma can be applied for the RRU model only when $m_1=m_2$.
For the case $m_1\neq m_2$ in an RRU model,
May and Flournoy (2009) established in (see~\cite{Flournoy.et.al.09}) that $\frac{Y_n}{n}\stackrel{a.s.}{\rightarrow}\max\{m_1;m_2\}$.
In the case $m_1=m_2$, we are able to establish that the limiting proportion $Z_{\infty}$
has no point mass within the open interval $(\rho_2,\rho_1)$.
This is stated in the following lemma.
\begin{lem}\label{lem:no_atoms}
Assume~\eqref{eq:rho_convergence_as} with $\rho_1>\rho_2$ and $m_1=m_2=m$.
Then, for any $x\in(\rho_2,\rho_1)$, we have $\bm{P}(Z_{\infty}=x)=0$.
\end{lem}
Point masses of probability are possible at values $\rho_1$ and $\rho_2$.

\subsection{Asymptotic distribution of the sampled balls}   \label{subsection_CLT}

The second order asymptotic results of the proportion of sampled balls are concerned with the concept of stable convergence (see~\cite{Hall.et.al.80}),
which provides a particularly elegant approach to martingale central limit theory.
Formally, let $\{\mathcal{X}_n;n\geq1\}$ be a random sequence on a probability space $(\Omega,\mathcal{F},\bm{P})$;
thus, we say that $\mathcal{X}_n\stackrel{d}{\rightarrow}\mathcal{X}$ (stably) if,
for every point $x$ of continuity for the cumulative distribution function of $\mathcal{X}$ and for every event $E\in\mathcal{F}$,
$$\lim_{n\rightarrow\infty}\ \bm{P}\left(\ \mathcal{X}_n\leq x,E\ \right)\ =\ \bm{P}\left(\ \mathcal{X}\leq x,E\ \right).$$

We now present the asymptotic distribution for the proportion of sampled balls in an RRU model.
Let us denote by $N_{1n}:=\sum_{i=1}^nX_i$ and $N_{2n}:=\sum_{i=1}^n(1-X_i)=n-N_{1n}$
the number of red and white balls, respectively, sampled form the urn up to time $n$.
Moreover, let $\sigma_1^2:=\bm{Var}[D_{1,1}]$ and $\sigma_2^2:=\bm{Var}[D_{2,1}]$.
The result is the following
\begin{thm}\label{thm:CLT_RRU}
Consider an RRU model and assume $m_1=m_2=m$. Then,
\begin{equation*}
\sqrt{n}\left(\frac{N_{1n}}{n}-Z_{\infty}\right)\ \stackrel{d}{\rightarrow}\ \mathcal{N}(0,\Sigma),\ \ \ \ (stably)
\end{equation*}
where
\begin{equation}\label{eq:CLT_RRU_sigma_2}
\Sigma\ :=\ \left(1+\frac{2\bar{\Sigma}}{m^2}\right)Z_{\infty}(1-Z_{\infty}),\ \qquad\ \bar{\Sigma}\ :=\ (1-Z_{\infty})\sigma^2_1+
Z_{\infty}\sigma^2_2.
\end{equation}
\end{thm}

We now present the asymptotic distribution for the proportion of sampled balls in an ARRU model.
This result can be derived by Theorem~\ref{thm:CLT_RRU} on the set of trajectories that do not cross the thresholds
$\hat{\rho}_{1,n}$ and $\hat{\rho}_{2,n}$ infinitely often,
and hence $\{Z_{\infty}\neq\{\rho_2,\rho_1\}\}$.
To this end, we introduce a sequence of random sets $\{A_n;n\geq1\}$ such that $A_n\in\mathcal{F}_n$ and
$A_n\subset A_{n+1}$ for any $n\geq1$, and $\cup_{n\geq1}A_n=(\rho_2,\rho_1)$.
In particular, we fix $0<\alpha< 1/2$ and we define $A_n$ as follows:
\begin{equation}\label{def:A_n_statement}
A_n\ :=\ \left(\ \rho_2+CY_n^{-\alpha}\ ,\ \rho_1-CY_n^{-\alpha}\ \right),
\end{equation}
where $0<C<\infty$ is a positive constant. The choice of $\{A_n;n\geq1\}$ in~\eqref{def:A_n_statement}
allows us to apply the estimates of Lemma~\ref{lem:same_mean} in the proof of the limit distribution, in order to obtain
the equivalence: $\{Z_n\in A_n,ev.\}=\{Z_{\infty}\in(\rho_2,\rho_1)\}$ a.s.,
where \textit{ev.} stands for \textit{eventually}, which means \textit{for all but a finite number of terms}.
The limit distribution for the ARRU model is expressed in the following result.
\begin{thm}\label{thm:CLT_ARRU}
Assume~\eqref{eq:rho_convergence_as} with $\rho_1>\rho_2$ and $m_1=m_2=m$. Then,
\begin{equation*}
\underline{\lim}_n \{Z_n\in A_n\}\ =\ \overline{\lim}_n \{Z_n\in A_n\}\ =\ \{Z_{\infty}\in (\rho_2,\rho_1)\},
\end{equation*}
and, on the sequence of sets $(\{Z_n\in A_n\}, n\geq1)$, we have
\begin{equation*}
\sqrt{n}\left(\frac{N_{1n}}{n}-Z_{\infty}\right)\ \stackrel{d}{\rightarrow}\ \mathcal{N}(0,\Sigma),\ \ \ \ (stably)
\end{equation*}
where, as in~\eqref{eq:CLT_RRU_sigma_2},
\begin{equation*}
\Sigma\ :=\ \left(1+\frac{2\bar{\Sigma}}{m^2}\right)Z_{\infty}(1-Z_{\infty}),\ \qquad\ \bar{\Sigma}\ :=\ (1-Z_{\infty})\sigma^2_1+
Z_{\infty}\sigma^2_2.
\end{equation*}
\end{thm}
It is worth noticing that the limiting distribution obtained in Theorem~\ref{thm:CLT_RRU} and Theorem~\ref{thm:CLT_ARRU} is not Gaussian but a mixture distribution.

As a corollary of the methods of proof of Theorem~\ref{thm:CLT_RRU} and Theorem~\ref{thm:CLT_ARRU} one can obtain the asymptotic distribution of
$\sqrt{n} (Z_n-Z_{\infty})$. We state this result without proof.
\begin{thm}\label{thm:CLT_ARRU_Z}
Assume~\eqref{eq:rho_convergence_as} with $\rho_1>\rho_2$ and $m_1=m_2=m$.
Then, conditionally on $\mathcal{F}_n$, on the sequence of sets $(\{Z_n\in A_n\}, n\geq1)$, we have
\begin{equation*}
\sqrt{n}\left(Z_n-Z_{\infty}\right)\ \stackrel{d}{\rightarrow}\ \mathcal{N}(0,\Sigma_Z),\ \ \ \ (stably)
\end{equation*}
where
\begin{equation*}
\Sigma_Z\ :=\ \left(1+\frac{\bar{\Sigma}}{m^2}\right)Z_{\infty}(1-Z_{\infty}),\ \qquad\ \bar{\Sigma}\ :=\ (1-Z_{\infty})\sigma^2_1+
Z_{\infty}\sigma^2_2.
\end{equation*}
\end{thm}

\section{Preliminary results}   \label{section_preliminary_results}

In this section, we present some preliminary estimates that are required to understand the dynamics of the ARRU model
and to prove the main results of the paper.
Most of the proofs of the results gathered by the literature are omitted, since the original proofs hold for all values of $m_1$ and $m_2$.

Initially, we show a useful expression of the excepted increments $(Z_{n+1}-Z_{n})$
conditionally to the story of the process $\mathcal{F}_n$, which is required to prove the consistency result
and in particular in the proof of Theorem~\ref{thm:increments_Z} in Section~\ref{section_WLLN}.
\begin{lem}\label{lem:mul_pag_sec}
For any $n\geq0$,
$$\bm{E}\left[Z_{n+1}-Z_{n}|\mathcal{F}_{n}\right]\ =\ Z_n(1-Z_n)B_n,$$
where
\begin{equation}\label{def:B_n}
B_n\ :=\ \bm{E}\left[\frac{D_{1,n+1}W_{1,n}}{Y_n+D_{1,n+1}W_{1,n}}-\frac{D_{2,n+1}W_{2,n}}{Y_n+D_{2,n+1}W_{2,n}}|\mathcal{F}_n\right].
\end{equation}
\end{lem}

\proof
The proof of this Lemma is based on a modification of the proof of Theorem 2 in~\cite{Muliere.et.al.06}.
First, note that, by definition
$$Z_{n+1}=X_{n+1}\frac{Y_{1,n}+D_{1,n+1}W_{1,n}}{Y_n+D_{1,n+1}W_{1,n}}+
(1-X_{n+1})\frac{Y_{1,n}}{Y_n+D_{2,n+1}W_{2,n}}$$
and since $X_{n+1}$ is conditionally to $\mathcal{F}_n$ independent of $D_{1,n+1}$ and $D_{2,n+1}$, we can get that
\[\begin{aligned}
\bm{E}[Z_{n+1}|\mathcal{F}_n]\ &&=&\ \bm{E}\left[Z_n\frac{Y_{1,n}+D_{1,n+1}W_{1,n}}{Y_n+D_{1,n+1}W_{1,n}}+
(1-Z_n)\frac{Y_{1,n}}{Y_n+D_{2,n+1}W_{2,n}}|\mathcal{F}_n\right]\\
&&=&\ \bm{E}\left[Z_n\left(\frac{Y_{1,n}+D_{1,n+1}W_{1,n}}{Y_n+D_{1,n+1}W_{1,n}}+
\frac{Y_{2,n}}{Y_n+D_{2,n+1}W_{2,n}}\right)|\mathcal{F}_n\right]
\end{aligned}\]
Analogously, we have that
$$\bm{E}[1-Z_{n+1}|\mathcal{F}_n]\ =\ \left[(1-Z_n)\left(\frac{Y_{2,n}+D_{2,n+1}W_{2,n}}{Y_n+D_{2,n+1}W_{2,n}}+
\frac{Y_{1,n}}{Y_n+D_{1,n+1}W_{1,n}}\right)|\mathcal{F}_n\right].$$
Therefore,
\[\begin{aligned}
&&&\bm{E}[Z_{n+1}-Z_n|\mathcal{F}_n]\ =\ \bm{E}[(1-Z_n)Z_{n+1}-Z_n(1-Z_{n+1})|\mathcal{F}_n]\\
&&=&\ Z_n(1-Z_n)\bm{E}\left[\frac{Y_{1,n}+D_{1,n+1}W_{1,n}}{Y_n+D_{1,n+1}W_{1,n}}+
\frac{Y_{2,n}}{Y_n+D_{2,n+1}W_{2,n}}\right.\\
&&&\left.-\frac{Y_{2,n}+D_{2,n+1}W_{2,n}}{Y_n+D_{2,n+1}W_{2,n}}-
\frac{Y_{1,n}}{Y_n+D_{1,n+1}W_{1,n}}|\mathcal{F}_n\right]\\
&&=&\ Z_n(1-Z_n)\bm{E}\left[\frac{D_{1,n+1}W_{1,n}}{Y_n+D_{1,n+1}W_{1,n}}-
\frac{D_{2,n+1}W_{2,n}}{Y_n+D_{2,n+1}W_{2,n}}|\mathcal{F}_n\right].
\end{aligned}\]
This concludes the proof.
\endproof

Now, we show that the number of balls sampled from the urn $N_{1,n}$, $N_{2,n}$ and the total number of balls in the urn $Y_n$, increase to infinity almost surely. To do that, we first need to show a lower bound for the increments of the process $Y_n$, which is given by the following:
\begin{lem}~\cite[Lemma 4.1]{Ghiglietti.et.al.15}\label{lem:lower_bound_Yn}
For any $i\geq1$, we have that
\begin{equation*}
\bm{E}\left[Y_i-Y_{i-1}|\mathcal{F}_{i-1}\right]\ \geq\ a\cdot\left(\frac{\min\{y_{1,0};y_{2,0}\}}{y_0+\left(i-1\right)b}\right).
\end{equation*}
\end{lem}

Here, we present the lemma on the divergence of the sequences $Y_n$, $N_{1,n}$ and $N_{2,n}$.
This result is obtained by using the conditional Borel-Cantelli lemma.
\begin{lem}~\cite[Lemma 4.2]{Ghiglietti.et.al.15}\label{lem:urn_not_stop}
Consider the urn model presented in Section~\ref{section_the_model}. Then,
\begin{itemize}
\item[(a)] $Y_n\stackrel{a.s.}{\rightarrow}\infty$;
\item[(b)] $\min\{N_{1,n};N_{2,n}\}\stackrel{a.s.}{\rightarrow}\infty$.
\end{itemize}
\end{lem}

The following lemma is needed in the proof of Theorem~\ref{thm:convergence_as}.
This result provides multiple equivalent ways to show the almost sure convergence of a real-valued process.
We consider a general real-valued process $\{Z_n;n\geq 0\}$ and two real numbers $d$ (down) and $u$ (up), with $d<u$.
The result requires two sequences of times $t_j(d,u)$ and $\tau_j(d,u)$ defined as follows:
for each $j\geq0$, $t_j(d,u)$ represents the time of the first up-cross of $u$ after $\tau_{j-1}(d,u)$, and
$\tau_j(d,u)$ represents the time of the first down-cross of $d$ after $t_j$.
Note that $t_j(d,u)$ and $\tau_j(d,u)$ are stopping times,
since the events $\{t_j(d,u)=k\}$ and $\{\tau_j(d,u)=k\}$ depend on $\{Z_n;n\leq k\}$,
which are measurable with respect to $\mathcal{F}_k$.
\begin{lem}~\cite[Theorem 2.1]{Aletti.et.al.13}\label{lem:Zn_converge_ifonlyif}
Let $\{Z_n;n\geq 0\}$ be a real-valued process in $\left[0,1\right]$.
Let $\tau_{-1}(d,u)=-1$ and define for every $j\geq0$ two stopping times
\begin{equation}\label{eq:def_tau}
\begin{aligned}
t_j(d,u) & =
\begin{cases}
\inf\{n>\tau_{j-1}(d,u):Z_n>u\} & \text{if }\{n>\tau_{j}(d,u):Z_n>u\} \neq\emptyset;
\\
+ \infty & \text{otherwise}.
\end{cases}
\\
\tau_j(d,u) & =
\begin{cases}
\inf\{n>t_{j}(d,u) :Z_n<d\}\ \  & \text{if }\{n>t_{j-1}(d,u) :Z_n<d\} \neq\emptyset;
\\
+ \infty & \text{otherwise}.
\end{cases}
\end{aligned}
\end{equation}
Then, the following three events are a.s. equivalent
\begin{itemize}
\item[(a)] $Z_n$ converges a.s.;
\item[(b)] for any $0<d<u<1$,
$$\lim_{j\rightarrow\infty} \bm{P}\left( t_{j}(d,u) < \infty \right) =0;$$
\item[(c)] for any $0<d<u<1$,
$$\sum_{j\geq1} \bm{P}\left( t_{j+1}(d,u) = \infty | t_j(d,u) < \infty\right) =\infty;$$
\end{itemize}
using the convention that $\bm{P}\left(t_{j+1}(d,u) = \infty  | t_j(d,u) < \infty\right) = 1$ when
$\bm{P}\left( t_j(d,u) = \infty\right)=1$.
\end{lem}

The following lemma provides lower bounds for the total number of balls in the urn at the times of up-crossings, $Y_{t_{j}}$.
The lemma gets used in the proof of Theorem~\ref{thm:convergence_as},
where conditioning to a fixed number of up-crossing ensures to have at least a number of balls $Y_n$ determined by the lower bounds of this lemma.
This result has been taken by Lemma 2.1 of~\cite{Aletti.et.al.13} and the proof is omitted since
the adaptive thresholds and the values of $m_1$ and $m_2$ do not play any role during up-crossings.
Hence, the proof reported in Lemma 2.1 of~\cite{Aletti.et.al.13} carries over to our model, with $D_n$ replaced by $Y_n$.
\begin{lem}~\cite[Lemma 2.1]{Aletti.et.al.13}\label{lem:Y_geometric_increasing}
For any $0<d<u<1$, we have that
$$Y_{t_{j}(d,u)}\geq \left(\frac{u\left(1-d\right)}{d\left(1-u\right)}\right)Y_{t_{j-1}(d,u)} \geq ... \geq \left(\frac{u\left(1-d\right)}{d\left(1-u\right)}\right)^{j}Y_{t_{0}(d,u)}.$$
\end{lem}

The following lemma provides a uniform bound for the generalized P\'{olya} urn with same reinforcement means,
which is needed in the proof of Theorem~\ref{thm:convergence_as}.
\begin{lem}~\cite[Lemma 3.2]{Aletti.et.al.13}\label{lem:same_mean}
Consider an RRU with $m_1 = m_2$. If $Y_0\geq 2b$,
then
$$\bm{P}\left( \sup_{n\geq1}|Z_n-Z_0| \geq h \right) \leq
\frac{b}{Y_0}\left(\frac{4}{h^2}+\frac{2}{h}\right)$$
for every $h>0$.
\end{lem}

Finally, we present an auxiliary result that provides an upper bound on the increments of the urn process $Z_n$,
by imposing a condition on the total number of balls in the urn $Y_n$.
\begin{lem}~\cite[Lemma 3.1]{Ghiglietti.et.al.15}\label{lem:Y_increments}
For any $\epsilon\in\left(0,1\right)$, we have that
\begin{equation}\label{eq:how_big_is_D}
\left\{\ Y_n > b \left(\frac{1-\epsilon}{\epsilon}\right)\ \right\}\ \ \ \subseteq\ \ \
\left\{\ |Z_{n+1}-Z_n|<\epsilon\ \right\}.
\end{equation}
\end{lem}

\section{Proof of weak consistency and related results}   \label{section_WLLN}

In this section, we prove the weak consistency for the urn proportion of the ARRU model, which is established in Theorem~\ref{thm:convergence_probability}.
This proof requires some probabilistic results concerning the ARRU model, which have been gathered in different subsections.
The proof of the weak consistency based on these results is then provided in Subsection~\ref{subsection_proof_WLLN}.

Let us start by describing the general structure of the proof.
The weak consistency is proved by showing that the process $\{\Delta_n;n\geq1\}$, defined as
\begin{equation}\label{def:increments_abs_values_urn_proportion}
\Delta_n:=|\rho_1-Z_n|,\ \ \forall\ n\geq0,
\end{equation}
converges to zero in probability.
To prove this, we want to exploit the fact that, unless $\Delta_n$ is arbitrarily close to zero,
the conditional expected increments of $\Delta_n$ are negative.
This result is obtained in Subsection~\ref{subsection_bound_1_similar_MRRU} by studying the conditional expected increments of $Z_n$.
Hence, to show that $\Delta_n$ is asymptotically close to zero, we need to investigate the expected increments of the process $\{\Delta_n;n\geq1\}$.
Since the increments of $\Delta_n$ are at the same order of $Y_n^{-1}$,
we first determine how fast the total number of balls in the urn, $Y_n$, increase to infinity.
This is addressed in Theorem~\ref{thm:Y_bounded_and_diverge}, where we show that
the total number of balls in the ARRU model increases linearly with the number of extractions from the urn.
For this reason, the increments of $\Delta_n$ are of the order of $n^{-1}$;
hence, we consider differences of $\Delta_n$ evaluated at linearly increasing times, i.e.
$G(n,c):=(\Delta_{n+nc}-\Delta_n)$,
such that the $L_1$ bounds obtained for such differences do not vanish as $n$ goes to infinity.
More specifically, we provide a negative upper bound for the expected differences $G(n,c)$,
which is not negligible unless $\Delta_n$ is asymptotically close to zero.
Formally, for any $\delta>0$, we show that for some $0<C<\infty$
\begin{equation}\label{thm:exponential_bounds_similar_MRRU}
\bm{E}\left[G(n,s_{\delta})\right]\ \leq\ -C\bm{P}(Q(\delta,n))\ +\ o(1),
\end{equation}
where $0<s_{\delta}<\infty$ is an appropriate constant and $Q(\delta,n):=\{\Delta_n>\delta\}$.
To obtain~\eqref{thm:exponential_bounds_similar_MRRU}, we prove that the expected differences $G(n,s_{\delta})$ are:
(i) negative for moderate values of $\Delta_n$
(see Theorem~\ref{thm:exponential_bounds_similar_MRRU_giu});
(ii) negligible for small values of $\Delta_n$
(see Theorem~\ref{thm:exponential_bounds_similar_MRRU_su}).
These results are derived using comparison arguments with specific auxiliary urn models.
Finally, in Subsection~\ref{subsection_proof_WLLN} we use~\eqref{thm:exponential_bounds_similar_MRRU} and other preliminary results
to establish the weak consistency.

\subsection{Harmonic moments of $Y_n$}   \label{subsection_harmonic_moments}

In this subsection, we establish that the total number of balls in the ARRU model increases linearly with
the number of extractions from the urn.
Moreover, this result ensures uniform bounds for the harmonic moments of the total number of balls.

Before presenting the main result, we introduce some notation.
For any $0<c\leq C<\infty$ and for all $n\geq0$, let $F_n(c,C)\in \mathcal{F}_n$ be the set defined as follows
\begin{equation*}
F_n(c,C):=\{y_0+cn\leq Y_n \leq y_0+Cn\}.
\end{equation*}
Here, we show that, for some $c$ and $C$, $\bm{P}(F_n(c,C))$ converges to one exponentially fast, which implies $\bm{P}(F^c_n(c,C),i.o.)=0$.
Moreover, this result provides uniform bounds for the moments of ${n}/{Y_n}$.
The following theorem makes this result precise.

\begin{thm}\label{thm:Y_bounded_and_diverge}
Under assumption~\eqref{ass:range_bounded}, for any $0<z_{\min}<\rho_{\min}$ and $\rho_{\max}<z_{\max}<1$,
there exists $\epsilon_{z}>0$ such that
\begin{equation}\label{eq:Z_exp_decay}
\bm{P}\left(z_{\min}\leq Z_n \leq z_{\max}\right)\ \geq\ 1-\exp(-\epsilon_{z}n).
\end{equation}
Moreover, there exist $0<c_1<C_1<\infty$ and $\epsilon_{y}>0$ such that
\begin{equation}\label{eq:Y_exp_decay}
\bm{P}(y_0+c_1n\leq Y_n \leq y_0+C_1n)\ \geq\ 1-\exp(-\epsilon_{y}n),
\end{equation}
for large $n$.
As a consequence, for any $j\geq1$
\begin{equation}\label{eq:Y_bounded}
\sup_{n\geq0}\left\{\bm{E}\left[\left(\frac{n}{Y_n}\right)^j\right]\right\}\ <\ \infty.
\end{equation}
\end{thm}

To ease notation in the rest of paper, we will refer to $F_n$ as
\begin{equation}\label{def:F_n}
F_n:=\{y_0+c_1n\leq Y_n \leq y_0+C_1n\},
\end{equation}
where $0<c_1<C_1<\infty$ are the constants determined in Theorem~\ref{thm:Y_bounded_and_diverge} to obtain~\eqref{eq:Y_exp_decay}.

\proof
Let $c_{\min}:=\min\{\rho_{\min};1-\rho_{\max}\}$, fix an arbitrary $0<c<c_{\min}$ and consider the following sets
\[\begin{aligned}
A_{d,n}\ &&:=&\ \left\{\bigcup_{n/2\leq i\leq n}\left\{Z_i<c\right\}\right\},\\
A_{c,n}\ &&:=&\ \left\{\bigcap_{n/2\leq i\leq n}\left\{c<Z_i<1-c\right\}\right\},\\
A_{u,n}\ &&:=&\ \left\{\bigcup_{n/2\leq i\leq n}\left\{Z_i>1-c\right\}\right\}.
\end{aligned}\]
In the proof of Theorem 3.1 in~\cite{Ghiglietti.et.al.15}, it is proved that $\bm{P}(A_{d,n})$ and $\bm{P}(A_{u,n})$
converges exponentially fast to zero, provided that for some $\epsilon_1,\epsilon_2>0$
\begin{equation}\label{ass:old_proof_harmonic_1}
\bm{P}\left(\hat{\rho}_{1,n}>\rho_1+\epsilon_1\right)\ \leq\ c_0\exp\left(-n\epsilon_1^2\right),
\end{equation}
and
\begin{equation}\label{ass:old_proof_harmonic_1}
\bm{P}\left(\hat{\rho}_{2,n}<\rho_2-\epsilon_2\right)\ \leq\ c_0\exp\left(-n\epsilon_2^2\right).
\end{equation}
Thus, setting $\epsilon_1$ and $\epsilon_2$ such that $\rho_1+\epsilon_1>\rho_{\max}$ and $\rho_2-\epsilon_2<\rho_{\min}$ and
using~\eqref{ass:range_bounded}, we can follow the same arguments obtaining that, for any $0<c<c_{\min}$,
$\bm{P}(A_{d,n})$ and $\bm{P}(A_{u,n})$ converges exponentially fast to zero,
which naturally implies~\eqref{eq:Z_exp_decay} since $z_{\min}<\rho_{\min}\leq c_{\min}$ and $z_{\max}>\rho_{\max}\geq 1-c_{\min}$.

Now, we prove~\eqref{eq:Y_exp_decay}. Since the reinforcements are a.s. bounded, i.e. $|D_{j,n}|<b$ for any $n\geq1$ and $j=1,2$,
we trivially have that $\bm{P}\left(Y_n\geq y_0 + nb\right)=0$. Thus, we will show the exponential decay of $\bm{P}\left(Y_n-y_0\leq c_1n\right)$.
Moreover, since from~\eqref{eq:Z_exp_decay} for any $0<c<c_{\min}$ there exists $\epsilon_{z}$ such that $\bm{P}\left(A_{c,n}\right)\geq1-\exp(-\epsilon_{z}n)$ ,
we will focus on the probability
$\bm{P}\left(\left\{Y_n-y_0\leq c_1n\right\}\cap\{A_{c,n}\}\right)$.

First, consider the following relation on the increments of the total number of balls
$$Y_i-Y_{i-1}=D_{1,i}X_iW_{1,i-1}+D_{2,i}\left(1-X_i\right)W_{2,i-1}\geq a\left[X_iW_{1,i-1}+\left(1-X_i\right)W_{2,i-1}\right]$$
Then, note that, on the set $A_{c,n}$, the random variables
$$X_iW_{1,i-1}+\left(1-X_i\right)W_{2,i-1},\ \ \ i=n/2,..,n$$
are, conditionally to the $\sigma$-algebra $\mathcal{F}_{i-1}$, Bernoulli with parameter with parameter greater than or equal to $c$.
Hence, if we introduce $\{B_i;i\geq1\}$ a sequence of i.i.d. Bernoulli random variable with parameter $c$,
\[\begin{aligned}
\bm{P}\left(\left\{Y_n-y_0\leq c_1n\right\}\cap\{A_{c,n}\}\right)\ &&\leq&\
\bm{P}\left(\left\{Y_n-Y_{n/2}\leq c_1n\right\}\cap\{A_{c,n}\}\right)\\
&&\leq&\ \bm{P}\left(\left\{a\sum_{i=n/2}^nB_i\leq c_1n\right\}\cap\{A_{c,n}\}\right)\\
&&\leq&\ \bm{P}\left(\left\{\sum_{i=n/2}^nB_i\leq \frac{c_1}{a}n\right\}\right).\\
\end{aligned}\]
Now, we want to use the Chernoff's bound for i.i.d. random variables in $[0,1]$ (see~\cite{Dembo.et.al.98}):
\begin{equation}\label{eq:bernoulli_chernoff_result}
\bm{P}\left( S_n\leq c_0\cdot\bm{E}[S_n]\right)\ \leq\ \exp\left(-\frac{(1-c_0)^2}{2}\cdot\bm{E}[S_n]\right),
\end{equation}
where $c_0\in (0,1)$ and $S_n=\sum_{i=n/2}^nB_i$.
In our case, we have $\bm{E}[S_n]=nc/2$ and so $c_0=2c_1/(ac)$. Hence, by choosing $c_1$ small enough we can obtain $c_0<1$ which let us apply Chernoff's bound.
This implies~\eqref{eq:Y_exp_decay}.

Finally, we get the harmonic moments as follows
\[\begin{aligned}\bm{E}\left[\left(\frac{n}{Y_n}\right)^j\right]\ &&=&\ \bm{E}\left[\left(\frac{n}{Y_n}\right)^j\ind_{F_n(c_1,C_1)}\right]+
\bm{E}\left[\left(\frac{n}{Y_n}\right)^j\ind_{F^c_n(c_1,C_1)}\right]\\
&&\leq&\ \bm{E}\left[\left(\frac{n}{y_0+c_1n}\right)^j\ind_{F_n(c_1,C_1)}\right]+
\left(\frac{n}{y_0}\right)^j\bm{E}\left[\ind_{F^c_n(c_1,C_1)}\right]\\
&&\leq&\ c_1^{-j}+y_0^{-j}n^j\exp(-\epsilon_{y}n).
\end{aligned}\]
\endproof

\subsection{$L_1$ Bound for the increments of $\Delta_n$}   \label{subsection_bound_1_similar_MRRU}

For any $\varepsilon>0$, let $R(\varepsilon,n):=\{|\hat{\rho}_{1,n}-\rho_1|<\varepsilon\}$
and $Q(\varepsilon,n):=\{\Delta_n>\varepsilon\}$, where
we recall from~\eqref{def:increments_abs_values_urn_proportion} that $\Delta_n=|\rho_1-Z_n|$.
The following result provides an upper bound on the increments of $\Delta_n$.

\begin{thm}\label{thm:increments_Z}
Let $m_1>m_2$ and assume~\eqref{ass:range_bounded} and~\eqref{eq:rho_convergence_probability}.
For any $\varepsilon>0$, there exists $0<c_2<\infty$ and a sequence of random variables $\{\psi_n;n\geq0\}$ with
$\bm{E}[|\psi_n|]=o(n^{-1})$, such that
\begin{equation}\label{eq:increments_Z}
\bm{E}\left[G(n,n^{-1})\ind_{Q(\varepsilon,n)}|\mathcal{F}_n\right]\ \leq
\ -n^{-1} \cdot c_2\ind_{Q(\varepsilon,n)}\ +\ \psi_n,
\end{equation}
where we recall $G(n,n^{-1})=\left(\Delta_{n+1}-\Delta_{n}\right)$.
\end{thm}

The behavior and the sign of the excepted increments of the urn proportion $G(n,n^{-1})$ required to prove Theorem~\ref{thm:increments_Z}
depend on the position of $Z_n$ respect to $\rho_1$.
For this reason, we study separately the cases when $Z_n$ is above or below $\rho_1$.
Formally, we define
\begin{equation}\label{eq:Q_minus_Q_plus}
Q^{-}(\varepsilon,n):=\{Z_n<\rho_1-\varepsilon\},\ \ \ \ Q^{+}(\varepsilon,n):=\{Z_n>\rho_1+\varepsilon\},
\end{equation}
so that $Q(\varepsilon,n)=Q^{+}(\varepsilon,n)\cup Q^{-}(\varepsilon,n)$.
Specifically, we present Lemma~\ref{lem:I_n_minus} and Lemma~\ref{lem:I_n_plus} that provide bounds for the expected increments $G(n,n^{-1})$
on the sets $Q^{-}(\varepsilon,n)$ and $Q^{+}(\varepsilon,n)$, respectively.
The proof of Theorem~\ref{thm:increments_Z} is presented after the proofs of Lemma~\ref{lem:I_n_minus} and Lemma~\ref{lem:I_n_plus}.
\begin{lem}\label{lem:I_n_minus}
Let $A_n\in\mathcal{F}_n$ be such that $A_n\subset Q^{-}(\varepsilon,n)$.
Then, we have that
\begin{equation}\label{eq:increments_Z_meno}
\bm{E}\left[\left(Z_{n+1}-Z_{n}\right)\ind_{A_n}\right]\ \geq\ n^{-1} \cdot c_2\bm{P}(A_{n})\ -\ o(n^{-1}).
\end{equation}
\end{lem}

\proof
Let $I_n:=\bm{E}\left[\left(Z_{n+1}-Z_{n}\right)\ind_{A_n}\right]$ and, since $A_n\in\mathcal{F}_n$, we can use Lemma~\ref{lem:mul_pag_sec} obtaining
\begin{equation}\label{eq:first_step_B_n}
I_n\ =\ \bm{E}\left[ \bm{E}\left[ Z_{n+1}-Z_n |\mathcal{F}_n\right]\ind_{A_n} \right]\ =\ \bm{E}\left[Z_n(1-Z_n)B_n\ind_{A_n}\right],
\end{equation}
where we recall that $B_n$ is defined in~\eqref{def:B_n} as follows
$$B_n\ :=\ \bm{E}\left[\frac{D_{1,n+1}W_{1,n}}{Y_n+D_{1,n+1}W_{1,n}}-\frac{D_{2,n+1}W_{2,n}}{Y_n+D_{2,n+1}W_{2,n}}\big|\mathcal{F}_n\right].$$
Now, note the following relation
$$\{Z_n\leq\hat{\rho}_{1,n}\}\supset Q^{-}(\varepsilon,n)\cap R(\varepsilon,n)$$
where $R(\varepsilon,n)=\{|\hat{\rho}_{1,n}-\rho_1|<\varepsilon\}$.
Since $A_n\subset Q^{-}(\varepsilon,n)$, on the set $A_n$ the previous relation becomes
$\{Z_n\leq\hat{\rho}_{1,n}\}\supset R(\varepsilon,n)$, which implies $W_{1,n}\geq \ind_{R(\varepsilon,n)}$.
Combining this argument with $W_{2,n}\leq 1$,
we obtain on the set $A_n$ the following inequality
$$B_n\ \geq\ \bm{E}\left[\left(\frac{D_{1,n+1}\ind_{R(\varepsilon,n)}}{Y_n+D_{1,n+1}\ind_{R(\varepsilon,n)}}-
\frac{D_{2,n+1}}{Y_n+D_{2,n+1}}\right)\big|\mathcal{F}_{n}\right].$$
Then, by using $D_{2,n+1}\geq0$ and $D_{1,n+1}\ind_{R(\varepsilon,n)}\leq b$ a.s., we obtain that, on the set $A_n$,
$$B_n\ \geq\ \bm{E}\left[\left(\frac{D_{1,n+1}\ind_{R(\varepsilon,n)}}{Y_n+b}-
\frac{D_{2,n+1}}{Y_n}\right)\big|\mathcal{F}_{n}\right]\ =\ E_{1n}\ -\ E_{2n},$$
where
$$E_{1n}:=\frac{m_1\ind_{R(\varepsilon,n)}-m_2}{Y_n+b},\ \mbox{and}\ E_{2n}:=\frac{m_2b}{Y_n(Y_n+b)}.$$
First, note that
$$\bm{E}\left[Z_n(1-Z_n)E_{2n}\ind_{A_n}\right]\ \leq\ \bm{E}\left[E_{2n}\right]\ \leq\ m_2b\sup_{k\geq1}\bm{E}\left[\left(\frac{k}{Y_k}\right)^2\right]
n^{-2}.$$
Now, using~\eqref{eq:Y_bounded} it follows that
$$\bm{E}\left[Z_n(1-Z_n)E_{2n}\ind_{A_n}\right]\ =\ O(n^{-2}).$$
Thus, from~\eqref{eq:first_step_B_n} we have
\begin{equation}\label{eq:second_step_B_n}
I_n\ \geq\ \bm{E}\left[Z_n(1-Z_n)E_{1n}\ind_{A_n}\right] - o(n^{-1}).
\end{equation}
Now, consider the set $F_n$ defined in~\eqref{def:F_n} as
$$F_n\ =\ \left\{\ c_1n\leq Y_n-y_0 \leq C_1n\ \right\},$$
where we recall that, by~\eqref{eq:Y_exp_decay} in Theorem~\ref{thm:Y_bounded_and_diverge}, $\bm{P}(F_n^c)\leq\exp\left(-\epsilon_{y}n\right)$.
Moreover, let $\ind_{A_n}=J_{1n}+J_{2n}$, where $J_{1n}:=\ind_{A_n\cap F_n}$ and $J_{2n}:=\ind_{A_n\cap F_n^c}$.
Thus, concerning $J_{2n}$ we have that
$$|\bm{E}\left[Z_n(1-Z_n)E_{1n}J_{2n}\right]|\leq \max_{n\geq0}\{|E_{1n}|\}\bm{P}\left( F_n^c\right)=o(n^{-1}),$$
since $\max_{n\geq0}\{|E_{1n}|\}\leq b/y_0$ a.s.
Thus, returning to~\eqref{eq:second_step_B_n} we have that
\begin{equation}\label{eq:third_step_B_n}
I_n\ \geq\ \bm{E}\left[Z_n(1-Z_n)E_{1n}J_{1n}\right] - o(n^{-1}).
\end{equation}
Now, consider the further decomposition $J_{1n}=J_{11n}+J_{12n}$,
where $J_{11n}:=\ind_{A_n\cap F_n\cap\{E_{1n}\geq0\}}$ and $J_{12n}:=\ind_{A_n\cap F_n\cap\{E_{1n}<0\}}$.
Thus, concerning $J_{12n}$ we have that
$$\bm{E}\left[Z_n(1-Z_n)E_{1n}J_{12n}\right]\ \geq\ -\left(\frac{m_2}{y_0+c_1(n+1)}\right)\bm{P}\left(A_n\cap\{E_{1n}<0\}\right);$$
moreover, since $\bm{P}(Z_n<z_{\min})$ and $\bm{P}(Z_n>z_{\max})$ converge to zero exponentially fast from~\eqref{eq:Z_exp_decay} in Theorem~\ref{thm:Y_bounded_and_diverge}, we obtain
$$\bm{E}\left[Z_n(1-Z_n)E_{1n}J_{11n}\right]\ \geq\ \left(\frac{z_{\min}\left(1-z_{\max}\right)(m_1-m_2)}{y_0+C_1(n+1)}\right)\bm{P}\left(A_n\cap\{E_{1n}>0\}\right)\ -\ o(n^{-1})$$
Therefore, from~\eqref{eq:third_step_B_n} we have
$$I_n\ \geq\ n^{-1}c_2\bm{P}\left(A_n\right)\ -\ O(n^{-1})\bm{P}\left(E_{1n}<0\right)\ -\ o(n^{-1}),$$
where $0<c_2<\infty$ is an appropriate constant.
Hence, since from $m_1>m_2$ we have $\{E_{1n}<0\}\equiv R^c(\varepsilon,n)$,
result~\eqref{eq:increments_Z_meno} is obtained by establishing $\bm{P}\left(E_{1n}<0\right)\rightarrow0$.
To this end, note that
$$\bm{P}\left(E_{1n}<0\right)\ =\ 1-\bm{P}\left(R(\varepsilon,n)\right)\ \rightarrow\ 0,$$
where $\bm{P}\left(R(\varepsilon,n)\right)\rightarrow1$ follows from $\hat{\rho}_1\stackrel{p}{\rightarrow}\rho_1$,
which is stated in~\eqref{eq:rho_convergence_probability} since $m_1>m_2$.
\endproof

Let us recall that from~\eqref{eq:Q_minus_Q_plus} $Q^{+}(\varepsilon,n)=\{Z_n>\rho_1+\varepsilon\}$.
We have the following result
\begin{lem}\label{lem:I_n_plus}
Let $A_n\in\mathcal{F}_n$ be such that $A_n\subset Q^{+}(\varepsilon,n)$.
Then, we have that
\begin{equation}\label{eq:increments_Z_piu}
\bm{E}\left[\left(Z_{n+1}-Z_{n}\right)\ind_{A_n}\right]\ \leq\ -n^{-1} \cdot c_2\bm{P}(A_{n})\ +\ o(n^{-1}).
\end{equation}
\end{lem}

\proof
The proof of this Lemma is obtained by following analogous arguments of the proof of Lemma~\ref{lem:I_n_minus}.
In fact, we can first apply Lemma~\ref{lem:mul_pag_sec}, then note that
$$\{Z_n\leq\hat{\rho}_{1,n}\}\subset Q^{+c}(\varepsilon,n)\cup R^c(\varepsilon,n),$$
and
$$\{Z_n\geq\hat{\rho}_{2,n}\}\supset Q^{+}(\varepsilon,n)\cap R(\varepsilon,n),$$
where we recall that $R(\varepsilon,n):=\{|\hat{\rho}_{1,n}-\rho_1|<\varepsilon\}$.
Hence, since $A_n\subset Q^{+}(\varepsilon,n)$, on the set $A_n$ we have that
$W_{1,n}\leq \ind_{R^c(\varepsilon,n)}$ and $W_{2,n}\geq \ind_{R(\varepsilon,n)}$, which lead to the following inequality
$$B_n\ \leq\ \bm{E}\left[\left(\frac{D_{1,n+1}\ind_{R^c(\varepsilon,n)}}{Y_n+D_{1,n+1}\ind_{R^c(\varepsilon,n)}}-
\frac{D_{2,n+1}\ind_{R(\varepsilon,n)}}{Y_n+D_{2,n+1}\ind_{R(\varepsilon,n)}}\right)|\mathcal{F}_{n}\right].$$
Then, by applying some standard calculations, we obtain that, on the set $A_n^+$,
\[\begin{aligned}
B_n\ &&\leq&\ \bm{E}\left[\left(\frac{D_{1,n+1}\ind_{R^c(\varepsilon,n)}}{Y_n}-
\frac{D_{2,n+1}\ind_{R(\varepsilon,n)}}{Y_n+b\ind_{R(\varepsilon,n)}}\right)|\mathcal{F}_{n}\right]\\
&&=&\ \frac{m_1\ind_{R^c(\varepsilon,n)}}{Y_n}-\frac{m_2\ind_{R(\varepsilon,n)}}{Y_n+b\ind_{R(\varepsilon,n)}}\\
&&=&\ \frac{m_1\ind_{R^c(\varepsilon,n)}-m_2\ind_{R(\varepsilon,n)}}{Y_n+b}.
\end{aligned}\]
Now, we can go through the same previous calculations using $\bm{P}\left( F_n^c\right)=o(n^{-1})$,~\eqref{thm:Y_bounded_and_diverge}
and $\bm{P}\left(R^c(\varepsilon,n)\right)\rightarrow0$, in order to prove~\eqref{eq:increments_Z_piu}.
\endproof

\proof[Theorem~\ref{thm:increments_Z}]
First, note that establishing~\eqref{eq:increments_Z} is equivalent to proving that for any $\mathcal{A}_n\in\mathcal{F}_n$ and letting
$A_n:=\mathcal{A}_n\cap Q_{(\varepsilon)n}$:
\begin{equation*}
\bm{E}\left[G(n,n^{-1})\ind_{A_n}\right]\ \leq\ -n^{-1} \cdot
c_2\bm{P}(A_{n})\ +\ o(n^{-1}),
\end{equation*}
where we recall that $G(n,n^{-1})=\left(\Delta_{n+1}-\Delta_n\right)$.
Hence, consider $A_{n}^+:=A_n\cap Q^{+}(\varepsilon,n)$ and $A_{n}^-:=A_n\cap Q^{-}(\varepsilon,n)$.
Since $A_{n}^+\cap A_{n}^-=\emptyset$ and $A_{n}^+\cup A_{n}^-=A_n$, we have the following decomposition
\begin{equation}\label{eq:decomposition_A_increments}
\bm{E}\left[G(n,n^{-1})\ind_{A_n}\right]\ =\ I_n^+\ -\ I_n^-,
\end{equation}
where
$$I_n^+\ :=\ \bm{E}\left[\left(Z_{n+1}-Z_n\right)\ind_{A_n^+}\right],\ \qquad\ I_n^-\ :=\ \bm{E}\left[\left(Z_{n+1}-Z_n\right)\ind_{A_n^-}\right].$$
By applying Lemma~\ref{lem:I_n_minus} and~\ref{lem:I_n_plus} to $I_n^-$ and $I_n^+$, respectively, we obtain
\[\begin{cases}
I_n^-\ \geq\ n^{-1} \cdot c_2\bm{P}(A_{n}^{-})\ -\ o(n^{-1}),\\
I_n^+\ \leq\ -n^{-1} \cdot c_2\bm{P}(A_{n}^{+})\ +\ o(n^{-1}).
\end{cases}\]
This concludes the proof.
\endproof

\subsection{$L_1$ Bound for $\Delta_n$ at linearly increasing times}   \label{subsection_bound_1_similar_MRRU_interval_time}

In this subsection, we provide an upper bound for the increments of $\Delta_n$
evaluated at linearly increasing times, i.e. $G(n,c)=(\Delta_{n+nc}-\Delta_n)$ and $c>0$,
where we recall from~\eqref{def:increments_abs_values_urn_proportion} that $\Delta_n=|\rho_1-Z_n|$.
To this end, we claim that, for any fixed $\delta>0$, there exist a value $c>0$ such that
$$\bm{P}\left(\ \left\{|Z_{n+ns_{\delta}}-Z_{n}|> \delta/2\right\}\ \cap\ F_n\ \right)\ =\ 0,$$
where we recall from~\eqref{def:F_n} that $F_n:=\{y_0+c_1n\leq Y_n \leq y_0+C_1n\}$.
We will denote by $s_{\delta}$ one of these values of $c$.

We can compute precisely the range of values admissible for $s_{\delta}$:
on the set $F_n$, we obtain
$$|Z_{n+nc}-Z_{n}|\leq b\sum_{i=n}^{n+nc}\frac{1}{Y_i}\leq
\frac{b}{c_1}\sum_{i=n}^{n+nc}\frac{1}{i}=\frac{b}{c_1}\log\left(1+c\right),$$
where we recall that $b$ is the maximum value of the urn reinforcements, i.e. $D_{1,n},D_{2,n}\leq b$ a.s. for any $n\geq1$.
Then, imposing $|Z_{n+nc}-Z_{n}|<\delta/2$, we obtain
\begin{equation}\label{def:s_delta}
s_{\delta}\ \in\ \left(\ 0\ ,\ \exp\left(\frac{c_1}{2b}\delta\right)-1\ \right).
\end{equation}
This ensures that $\bm{P}\left(\left\{|Z_{n+ns_{\delta}}-Z_{n}|> \delta/2\right\}\cap F_n\right)=0$.

The next theorem provides an $L_1$ upper bound for the difference $G(n,s_{\delta})=(\Delta_{n+ns_{\delta}}-\Delta_n)$
on the set $Q(\delta,n)=\{\Delta_n>\delta\}$.
An $L_1$ upper bound on the set $Q^c(\delta,n)$ is presented in Theorem~\ref{thm:exponential_bounds_similar_MRRU_su}.

\begin{thm}\label{thm:exponential_bounds_similar_MRRU_giu}
Let $m_1>m_2$,~\eqref{ass:range_bounded} and~\eqref{eq:rho_convergence_probability}.
Then, for any $\delta>0$ there exists a constant $0<C<\infty$ such that
\begin{equation}\label{eq:exponential_bounds_similar_MRRU_giu}
\bm{E}\left[G(n,s_{\delta})\ind_{Q(\delta,n)}\right]\ \leq\ -C\bm{P}(Q(\delta,n))\ +\ o(1).
\end{equation}
\end{thm}

\proof
First, note that using~\eqref{eq:Y_exp_decay} in Theorem~\ref{thm:Y_bounded_and_diverge}, we have
$$\left|\bm{E}\left[G(n,s_{\delta})\ind_{Q(\delta,n)\cap F_n^c}\right]\right|\ \leq\
\bm{P}(F_n^c)\ \rightarrow\ 0.$$
Hence, define
\begin{equation*}\label{def:G_n}
G_n\ :=\ \bm{E}\left[G(n,s_{\delta})\ind_{Q(\delta,n)\cap F_n}\right],
\end{equation*}
and consider the following expression
\begin{equation}\label{eq:G_n}
G_n\ =\ \sum_{i=n}^{n+ns_{\delta}-1}\bm{E}\left[G(i,i^{-1})\ind_{Q(\delta,n)\cap F_n}\right],
\end{equation}
where we recall that $G(i,i^{-1})=(\Delta_{i+1}-\Delta_i)$.
From the definition of $s_{\delta}$ in~\eqref{def:s_delta}, on the set $F_n$ we have that for all $i\in\{n,..,n+ns_{\delta}\}$
$$Q(\delta,n)\ \subset\ Q(\delta/2,i),$$
where we recall that $Q(\delta,n)=\{\Delta_n>\delta\}$ and $Q(\delta/2,i)=\{\Delta_i>\delta/2\}$.
Hence, by applying Theorem~\ref{thm:increments_Z} to each term of the sum in~\eqref{eq:G_n}, since $Q(\delta,n)\cap F_n\in\mathcal{F}_i$
for all $i\in\{n,..,n+ns_{\delta}\}$, we obtain
\[\begin{aligned}
\bm{E}\left[G(i,i^{-1})\ind_{Q(\delta,n)\cap F_n}\right]\ &&=&\
\bm{E}\left[\bm{E}\left[G(i,i^{-1})\ind_{Q(\delta/2,i)}|\mathcal{F}_i\right]\ind_{Q(\delta,n)\cap F_n}\right]\\
&&\leq&\ \bm{E}\left[\left(-i^{-1} \cdot c_2\ind_{Q(\delta/2,i)}\ +\ \psi_i\right)\ind_{Q(\delta,n)\cap F_n}\right]\\
&&=&\ -i^{-1} \cdot c_2\bm{P}\left(Q(\delta,n)\cap F_n\right)\ +\ \bm{E}\left[\psi_i\ind_{Q(\delta,n)\cap F_n}\right].
\end{aligned}\]
Now, note that from~\eqref{eq:Y_exp_decay} in Theorem~\ref{thm:Y_bounded_and_diverge} we have that
$\bm{P}\left(Q(\delta,n)\cap F_n\right)=\bm{P}\left(Q(\delta,n)\right)-o(i^{-1})$;
moreover, from Theorem~\ref{thm:increments_Z} $|\bm{E}\left[\psi_i\ind_{Q(\delta,n)\cap F_n}\right]|\leq \bm{E}\left[|\psi_i|\right]=o(i^{-1})$.
Thus, from~\eqref{eq:G_n} we have that
\[\begin{aligned}
G_n\ &&\leq&\ -\sum_{i=n}^{n+ns_{\delta}-1}i^{-1} \cdot c_2\bm{P}(Q(\delta,n))\ +\ \sum_{i=n}^{n+ns_{\delta}-1}o\left(i^{-1}\right)\\
&&\leq&\ -\log\left(1+s_{\delta}\right)\cdot c_2\bm{P}(Q(\delta,n))\ +\ o(1).
\end{aligned}\]
The result follows after calling $C:=c_2\log\left(1+s_{\delta}\right)$.
\endproof

Now, we show that the expected difference $G(n,s_{\delta})$
is asymptotically non-positive on the set $Q^c(\delta,n)$, for any $\delta>0$,
where we recall that $G(n,s_{\delta})=(\Delta_{n+ns_{\delta}}-\Delta_n)$, $Q(\delta,n)=\{\Delta_n>\delta\}$ and $\Delta_n=|\rho_1-Z_n|$.
The result is stated precisely in the following theorem.

\begin{thm}\label{thm:exponential_bounds_similar_MRRU_su}
Let $m_1>m_2$,~\eqref{ass:range_bounded} and~\eqref{eq:rho_convergence_probability}.
Then, for any $\delta>0$,
\begin{equation}\label{eq:exponential_bounds_similar_MRRU_su}
\overline{\lim}_n\bm{E}\left[G(n,s_{\delta})\ind_{Q^c(\delta,n)}\right]\ \leq\ 0.
\end{equation}
\end{thm}

To prove Theorem~\ref{thm:exponential_bounds_similar_MRRU_su}, we need to compare the ARRU model with two new urn models:
$\{\tilde{Z}^{+}_n;n\geq1\}$ and $\{\tilde{Z}^{-}_n;n\geq1\}$.
The dynamics of these processes is based on a sequence of random times
$\{t_n;n\geq1\}$ which describes relation between the process $\{\Delta_{n};n\geq1\}$ and an arbitrary fixed value $\nu>0$.
Specifically, fix $\nu>0$ and, for any $n\geq0$, define the set
$$\mathcal{T}_n\ :=\ \left\{0\leq k\leq n\ :\ Q^c(\nu,n-k)\right\},$$
where we recall $Q^c(\nu,n-k)=\{\Delta_{n-k}\leq\nu\}$.
Let $\{t_n;n\geq1\}$ be the sequence of random times defined as
\begin{equation}\label{def:t_n_lemma}
\begin{aligned}
t_n & =
\begin{cases}
\inf\{\mathcal{T}_n\} & \text{if }\mathcal{T}_n \neq\emptyset;\\
\infty & \text{otherwise}.
\end{cases}
\end{aligned}\end{equation}
The time $(n-t_n)$ indicates the last time up to $n$ the urn proportion is in the interval $(\rho_1-\nu,\rho_1+\nu)$.

First, let us describe the urn model $\{\tilde{Z}^{-}_n;n\geq1\}$.
Let $\tilde{I}^{-}=1$, $\tilde{y}_0\in(0,y_0)$ and $\tilde{z}_0^{-}\in(0,\rho_1-\nu)$.
The process $\{\tilde{Z}^{-}_n;n\geq1\}$, $\tilde{Z}^{-}_n=\tilde{Y}_{1,n}/(\tilde{Y}_{1,n}+\tilde{Y}_{2,n})$, evolves as follows:
if $t_{n-1}=0$, i.e. $\Delta_{n-1}\leq\nu$, or $t_{n-1}=\infty$, then $\tilde{X}_n=\ind_{\{U_{n}<\tilde{z}_0^{-}\}}$ and
\begin{equation}\label{def:urn_tilde_coupled_induction_1}\left\{
\begin{array}{l}
\tilde{Y}_{1,n}\ =\ \tilde{z}_0^{-}\cdot\tilde{y}_0\ +\ \tilde{X}_{n}D_{1,n}\tilde{I}^{-},\\
\\
\tilde{Y}_{2,n}\ =\ \left(1-\tilde{z}_0^{-}\right)\cdot \tilde{y}_0+\left(1-\tilde{X}_{n}\right)D_{2,n};
\end{array}
\right.\end{equation}
if $t_{n-1}=k\geq1$, i.e. $\Delta_{n-1}>\nu$, then $\tilde{X}_n=\ind_{\{U_{n}<\tilde{Z}_{n-1}\}}$ and
\begin{equation}\label{def:urn_tilde_coupled_induction_2}\left\{
\begin{array}{l}
\tilde{Y}_{1,n}\ =\ \tilde{Y}_{1,n-1}+\tilde{X}_{n}D_{1,n}\tilde{I}^{-},\\
\\
\tilde{Y}_{2,n}\ =\ \tilde{Y}_{2,n-1}+\left(1-\tilde{X}_{n}\right)D_{2,n}
\end{array}
\right.\end{equation}
then, $\tilde{Y}_{n}:=\tilde{Y}_{1,n}+\tilde{Y}_{2,n}$ and $\tilde{Z}_{n}:=\tilde{Y}_{1,n}/\tilde{Y}_{n}$.
The urn model is well defined since $t_{n-1}$ is $\mathcal{F}_{n-1}$-measurable.

Analogously, the urn model $\{\tilde{Z}^{+}_n;n\geq1\}$, $\tilde{Z}^{+}_n=\tilde{Y}_{1,n}/(\tilde{Y}_{1,n}+\tilde{Y}_{2,n})$, is defined by the same equations~\eqref{def:urn_tilde_coupled_induction_1} and~\eqref{def:urn_tilde_coupled_induction_2},
with $\tilde{I}^{-}$ and $\tilde{z}_0^{-}$ are replaced by $\tilde{I}^{+}=1$ and $\tilde{z}_0^{+}\in(\rho_1+\nu,1)$, respectively.

In the next lemma, we state an important relation among the processes $\{\tilde{Z}^{-}_n;n\geq1\}$, $\{\tilde{Z}^{+}_n;n\geq1\}$ and the urn proportion
of the ARRU model $\{Z_n;n\geq1\}$.
This result is needed in the proof of Theorem~\ref{thm:exponential_bounds_similar_MRRU_su}.
To ease calculations, let $h>0$ and fix the initial proportions $\tilde{z}_0^{-}$ and $\tilde{z}_0^{+}$ as follows:
\begin{equation}\label{eq:constraint_z_0}
\rho_1-\tilde{z}_0^{-}\ =\ \tilde{z}_0^{+}-\rho_1\ =\ \nu+h.
\end{equation}
Let $M_n:=\sum_{i=n}^{n+ns_{\delta}}\ind_{R^c(\nu,n)}$ and, for any $\epsilon>0$ define the set
\begin{equation}\label{def:M_n_epsilon}
\mathcal{M}_{n}^{\epsilon}\ :=\ \{M_n<ns_{\delta}\epsilon\},
\end{equation}
where we recall that $R(\nu,n)=\{|\hat{\rho}_{1,n}-\rho_1|\leq\nu\}$,
$s_{\delta}$ is such that
$\bm{P}\left(\left\{|G(n,s_{\delta})|> \delta/2\right\},F_n\right)=0$,
with $F_n=\{y_0+c_1n\leq Y_n \leq y_0+C_1n\}$ from~\eqref{def:F_n}.
Moreover, for any $n\geq1$ and $k\in\{n,..,n+ns_{\delta}\}$ let us define the set
\begin{equation}\label{def:E_nk}
E(n,k)\ :=\ \cup_{j=n}^{k}Q^c(\nu,j)\ \equiv\ \left\{\ \exists j\in\{n,..,k\}:\{\Delta_j\leq\nu\}\ \right\}.
\end{equation}
We also introduce the following notation:
$\tilde{\Delta}^{-}_{l}:=|\rho_1-\tilde{Z}^{-}_{l}|$, $\tilde{\Delta}^{+}_{l}:=|\rho_1-\tilde{Z}^{+}_{l}|$ and
$\tilde{\Delta}^{*}_{l}:=\max\left\{\tilde{\Delta}^{-}_{l},\tilde{\Delta}^{+}_{l}\right\}$.
Thus, we have the following result:

\begin{lem}\label{lem:induction2}
Let $m_1>m_2$,~\eqref{ass:range_bounded} and~\eqref{eq:rho_convergence_probability}.
Fix $n\geq1$, $\tilde{y}_0\in(0,y_0+c_1n)$, $\tilde{z}_0^{-}$ and $\tilde{z}_0^{+}$ as in~\eqref{eq:constraint_z_0}.
Consider the set $\mathcal{M}_{n}^{\epsilon}$ as defined in~\eqref{def:M_n_epsilon} with
\begin{equation}\label{eq:epsilon_induction_tilde}
0<\epsilon<\frac{c_1h}{bs_{\delta}}.
\end{equation}
Then, for any $n\geq 1$ and $l_n\in\{n+1,..,n+ns_{\delta}\}$, on the set $\mathcal{M}_{n}^{\epsilon}\cap F_n$ we have that
\begin{equation}\label{eq:induction2}
E(n,l_n) \subset\ Q^c(\tilde{\Delta}^{*}_l,l)\ a.s.,
\end{equation}
for all $l\in\{l_n+1,..,n+ns_{\delta}\}$.
\end{lem}

\proof
The proof will be by induction on $l\in\{l_n+1,..,n+ns_{\delta}\}$.
First, note that, from the definition of $\{t_n;n\geq1\}$ in~\eqref{def:t_n_lemma} and $E(n,k)$ in~\eqref{def:E_nk}, we always have
$$\{t_{l-1}=\infty\}\ \cap\ E(n,l_n)\ =\ \emptyset.$$
Hence, we never consider in this proof the set $\{t_{l-1}=\infty\}$.

Then, consider the set $\{t_{l-1}=0\}$ and note that, from the definition of $t_n$ in~\eqref{def:t_n_lemma}, $\{t_{l-1}=0\}\equiv Q^c(\nu,l-1)$,
which implies that, on the set $\{t_{l-1}=0\}\cap\{X_{l}=0\}$,
\begin{equation}\label{eq:hp_initial_induction_tilde_meno}
Z_{l}\ \geq\ \frac{(\rho_1-\nu) Y_{l-1}}{Y_{l-1}+D_{2,l}W_{2,l-1}}\ \geq\ \frac{\tilde{z}^{-}_0 \tilde{y}_0}{\tilde{y}_0+D_{2,l}}\ =\ \tilde{Z}_l^{-}\ \ \ \ a.s.,
\end{equation}
and, on the set $\{t_{l-1}=0\}\cap\{X_{l}=1\}$,
\begin{equation}\label{eq:hp_initial_induction_tilde_piu}
Z_{l}\ \leq\ \frac{(\rho_1+\nu) Y_{l-1}+D_{1,l}W_{1,l-1}}{Y_{l-1}+D_{1,l}W_{1,l-1}}\ \leq\ \frac{\tilde{z}^{+}_0 \tilde{y}_0+D_{1,l}}{\tilde{y}_0+D_{1,l}}\ =\ \tilde{Z}_l^{+}\ \ \ \ a.s.
\end{equation}
From~\eqref{eq:hp_initial_induction_tilde_meno} and~\eqref{eq:hp_initial_induction_tilde_piu} we have $\tilde{Z}_l^{-}\leq Z_{l}\leq\tilde{Z}_l^{+}$ a.s.,
that ensures that~\eqref{eq:induction2} is verified whenever $\{t_{l-1}=0\}$.

To prove~\eqref{eq:induction2} on the set $\{1\leq t_{l-1}<\infty\}$, we will show that,
defining $\tilde{A}_l^{-}:=\left\{\tilde{Z}_l^{-}\leq Z_{l}\right\}$, $\tilde{A}_l^{+}:=\left\{Z_{l}\leq\tilde{Z}_l^{+}\right\}$
and $B_n:=\mathcal{M}_{n}^{\epsilon}\cap F_n \cap \{1\leq t_{l-1}<\infty\}$,
\begin{equation}\label{eq:induction2_in_lemma}
\bm{P}\left(
\left\{\tilde{A}_l^{-}\cap Q^{-}(\nu,l-t_{l-1})\right\}\ \cup
\left\{\tilde{A}_l^{+}\cap Q^{+}(\nu,l-t_{l-1})\right\}\ |B_n\right)=1,
\end{equation}
for any $l\in\{l_n+1,..,n+ns_{\delta}\}$.
Moreover, from the definition of $\{t_n;n\geq1\}$ in~\eqref{def:t_n_lemma}, on the set $\{1\leq t_{l-1}<\infty\}$,
we note that
$$\{X_{l-t_{l-1}}=1\}\ \equiv\ \{Z_{l-t_{l-1}}\geq\rho_1+\nu\}=Q^{+}(\nu,l-t_{l-1}),$$
$$\{X_{l-t_{l-1}}=0\}\ \equiv\ \{Z_{l-t_{l-1}}\leq\rho_1-\nu\}=Q^{-}(\nu,l-t_{l-1}).$$
Hence, showing~\eqref{eq:induction2_in_lemma} is equivalent to establish the following
\begin{equation}\label{eq:induction2_in_lemma2}
\bm{P}\left(
\left\{\tilde{A}_l^{-}\cap \{X_{l-t_{l-1}}=0\}\right\}\ \cup
\left\{\tilde{A}_l^{+}\cap \{X_{l-t_{l-1}}=1\}\right\}\ |B_n\right)=1,
\end{equation}
Now, consider $\{1\leq t_{l-1}<\infty\}\cap\{X_{l-t_{l-1}}=0\}$, and by inductive hypothesis let $\omega$ belongs to the set
\begin{equation}\label{eq:hp_induction_tilde_meno}
\bigcap_{i=l-t_{l-1}}^{l-1}\tilde{A}_i^{-},
\end{equation}
where we recall that $\tilde{A}_i^{-}=\left\{\tilde{Z}^{-}_i\leq Z_i\right\}$.
Note that by~\eqref{eq:hp_initial_induction_tilde_meno} it follows that, on the set $\{t_{l-1}=1\}$,
condition~\eqref{eq:hp_induction_tilde_meno} is verified.
Hence, the result is achieved by establishing that~\eqref{eq:hp_induction_tilde_meno} implies $\omega$ belongs to $\tilde{A}_l^{-}$.

To this end, consider
$$Z_{l}=\frac{Z_{l-t_{l-1}-1}Y_{l-t_{l-1}-1}+\sum_{i=l-t_{l-1}+1}^{l}X_iD_{1,i}W_{1,i-1}}
{Y_{l-t_{l-1}-1}+\sum_{i=l-t_l+1}^{l}X_iD_{1,i}W_{1,i-1}+\sum_{i=l-t_{l-1}}^{l}\left(1-X_i\right)D_{2,i}W_{2,i-1}}.$$
Now, note that by~\eqref{eq:hp_induction_tilde_meno} we have
$X_{i}=\ind_{\{U_{i}<Z_{i-1}\}}\geq\ind_{\{U_{i}<\tilde{Z}_{i-1}^{-}\}}=\tilde{X}^{-}_{i+1}$ for any $i=l-t_{l-1}+1,...,l$.
Moreover, since $Z_{l-t_{l-1}-1}\geq\rho_1-\nu$, $Y_{l-t_{l-1}-1}\geq \tilde{y}_0$ and $X_{l-t_{l-1}}=0$
it follows that
$$Z_{l}\ \geq\ \frac{(\rho_1-\nu)\tilde{y}_0+\sum_{i=l-t_{l-1}+1}^{l}\tilde{X}^{-}_iD_{1,i}W_{1,i-1}}
{\tilde{y}_0+\sum_{i=l-t_l+1}^{l}\tilde{X}^{-}_iD_{1,i}W_{1,i-1}+\sum_{i=l-t_{l-1}}^{l}\left(1-\tilde{X}^{-}_i\right)D_{2,i}W_{2,i-1}}.$$
Note that, letting $n_0$ such that $\bm{P}(R(\nu,n_0))>\eta>0$, for any $n\geq n_0$ we have the following relation
$$\{Z_n\leq\hat{\rho}_{1,n}\}\supset Q^{-}(\nu,n)\cap R(\nu,n),$$
where we recall that $R(\nu,n)=\{|\hat{\rho}_{1,n}-\rho_1|<\nu\}$ and $Q^{-}(\nu,n)=\{Z_n<\rho_1-\nu\}$.
Hence, by definition of $t_{l-1}$ in~\eqref{def:t_n_lemma}, we have $Q^{-}(\nu,i)$ for any $i=l-t_{l-1},...,l-1$, and
$\{Z_i\leq\hat{\rho}_{1,i}\}\supset R(\nu,i)$, which implies $W_{1,i}\geq \ind_{R(\nu,i)}$.
Combining this argument with $W_{2,i}\leq 1$, we have that
$$Z_{l}\ \geq\ \frac{(\rho_1-\nu)\tilde{y}_0+\sum_{i=l-t_{l-1}+1}^{l}\tilde{X}^{-}_iD_{1,i}\ind_{R(\nu,i-1)}}
{\tilde{y}_0+\sum_{i=l-t_{l-1}+1}^{l}\tilde{X}^{-}_iD_{1,i}\ind_{R(\nu,i-1)}+\sum_{i=l-t_{l-1}+1}^{l}\left(1-\tilde{X}^{-}_i\right)D_{2,i}}.$$
In addition, on the set $\mathcal{M}_{n}^{\epsilon}$ we have that
\[\begin{aligned}\sum_{i=l-t_{l-1}+1}^{l}\tilde{X}^{-}_iD_{1,i}\ind_{R(\nu,i-1)}\ &&\geq&\
\sum_{i=l-t_{l-1}+1}^{l}\tilde{X}^{-}_iD_{1,i}-bM_n\\
&&\geq&\ \sum_{i=l-t_{l-1}+1}^{l}\tilde{X}^{-}_iD_{1,i}-nbs_{\delta}\epsilon.
\end{aligned}\]
Moreover, condition~\eqref{eq:epsilon_induction_tilde} ensures that
$$(\rho_1-\nu)\tilde{y}_0-nbs_{\delta}\epsilon\ \geq\ \tilde{z}^{-}_0\tilde{y}_0,$$
which implies $Z_{l}\geq\tilde{Z}^{-}_{l}$.
This concludes the proof of $\{\tilde{Z}^{-}_{l}\leq Z_{l}\}$.

Analogous arguments can be followed when we consider $\{1\leq t_{l-1}<\infty\}\cap\{X_{l-t_{l-1}}=1\}$.
In this case, by inductive hypothesis let $\omega$ belongs to the set
\begin{equation}\label{eq:hp_induction_tilde_piu}
\bigcap_{i=l-t_{l-1}}^{l-1}\tilde{A}_i^{+},
\end{equation}
where $\tilde{A}_i^{+}=\left\{\tilde{Z}^{+}_i\geq Z_i\right\}$.
Then, note that condition~\eqref{eq:hp_induction_tilde_piu} is verified for $t_{l-1}=1$ using~\eqref{eq:hp_initial_induction_tilde_piu}.
Hence, the result can be achieved by establishing in an analogous way that~\eqref{eq:hp_induction_tilde_piu} implies $\omega$ belongs to $\tilde{A}_l^{+}$.

Finally, combining $\tilde{A}_l^{-}$ and $\tilde{A}_l^{+}$, we obtain~\eqref{eq:induction2_in_lemma2}.
This concludes the proof.
\endproof

In the next lemma, we show an important result required in the proof of Theorem~\ref{thm:exponential_bounds_similar_MRRU_su},
concerning the probability that $\tilde{Z}_n$ exceeds an arbitrary threshold $l>0$.
This result is obtained by using comparison arguments between the process $\{\tilde{\Delta}^{*}_n;n\geq1\}$ and
the urn proportion of an RRU model, where we recall that $\tilde{\Delta}^{*}_n=\max\{\tilde{\Delta}^{-}_{l},\tilde{\Delta}^{+}_{l}\}$,
$\tilde{\Delta}^{-}_{l}:=|\rho_1-\tilde{Z}^{-}_{l}|$ and $\tilde{\Delta}^{+}_{l}:=|\rho_1-\tilde{Z}^{+}_{l}|$.
The result is the following,
\begin{lem}\label{lem:comparison_tilde_RRU}
Let $m_1>m_2$, and
\begin{equation}\label{def:Ttilde_n_H_n}
\tilde{T}_n\ :=\ \left\{k_n<t_n<\infty\right\},\ \qquad\ H_n\ :=\ \left\{\tilde{\Delta}^{*}_n>\nu\right\},
\end{equation}
where $\{k_n;n\geq1\}$ is a deterministic sequence such that $k_n\rightarrow\infty$.
Fix $0<\tilde{y}_0<\infty$ and define $\tilde{z}_0^{-}$ and $\tilde{z}_0^{+}$ as in~\eqref{eq:constraint_z_0}.
Then,
\begin{equation}\label{eq:result_comparison_tilde_RRU}
\lim_{n\rightarrow\infty}\bm{P}\left(H_n\cup\tilde{T}_n\right) =0.
\end{equation}
\end{lem}

\proof
Since $H_n=H_n^-\cup H_n^+$ where
$$H_n^-:=\left\{\tilde{Z}^-_n<\rho_1-\nu\right\},\ \mbox{and}\
H_n^+:=\left\{\tilde{Z}^+_n>\rho_1+\nu\right\},$$
equation~\eqref{eq:result_comparison_tilde_RRU} is established by proving
$$\lim_{n\rightarrow\infty}\bm{P}\left(H^-_n\cup\tilde{T}_n\right)+\bm{P}\left(H^+_n\cup\tilde{T}_n\right) =0.$$
We will show that $\bm{P}\left(H^-_n\cup\tilde{T}_n\right)\rightarrow0$,
since the proof of $\bm{P}\left(H^-_n\cup\tilde{T}_n\right)\rightarrow0$ is analogous.

First, we recall that $t_n$, defined in~\eqref{def:t_n_lemma}, satisfies that $Q^c(\nu,n-t_n)=\{\Delta_{n-t_n}\leq\nu\}$ and
when $t_n>0$, $Q(\nu,i)=\{\Delta_i>\nu\}$ for any $n-t_n< i\leq n$.
Hence, on the set $\tilde{T}_n$ the process $\tilde{Z}^-_i$ evolves at times $n-t_n< i\leq n$ as described in~\eqref{def:urn_tilde_coupled_induction_2},
yielding $\tilde{X}_{i}=\ind_{\{U_i<\tilde{Z}_{i-1}^-\}}$ and
\begin{equation}\label{eq:urn_tilde_coupled}\left\{
\begin{array}{l}
\tilde{Y}^-_{1,n}\ =\ \tilde{z}^-_{0}\tilde{y}_0\ +\ \sum_{i=n-t_n+1}^{n}\tilde{X}_{i}D_{1,i},\\
\\
\tilde{Y}^-_{2,n}\ =\ (1-\tilde{z}^-_{0})\tilde{y}_0\ +\ \sum_{i=n-t_n+1}^{n}\left(1-\tilde{X}_{i}\right)D_{2,i}.
\end{array}
\right.\end{equation}

Now, consider an RRU model $\{Z^{R}_{j};j\geq1\}$ with initial composition $(\tilde{y}^R_{1,0},\tilde{y}^R_{2,0})=(\tilde{z}^-_{0}\tilde{y}_0,(1-\tilde{z}^-_{0})\tilde{y}_0)$;
the reinforcements are defined as $D^{R}_{1,j}=D_{1,n-t_n+j}$ and $D^{R}_{2,j}=D_{2,n-t_n+j}$ for any $i\geq1$ a.s.;
the sampling process is modeled by $X^R_{j}:=\ind_{\{U^R_{j}<Z^R_{j-1}\}}$ and $U^{R}_{j}=U_{n-t_n+j}$ a.s.,
Hence, the composition of the RRU model at time $j\geq1$ can be expressed as follows:
\begin{equation}\label{eq:urn_RRU_coupled}
\begin{aligned}
Y^R_{1,j}\ &&=&\ \tilde{y}^R_{1,0}\ +\ \sum_{i=1}^{j}X_{n-t_n+i}D_{1,n-t_n+i}\\
&&=&\ \tilde{z}^-_{0}\tilde{y}^-_0\ +\ \sum_{i=n-t_n+1}^{n-t_n+j}X_{i}D_{1,i},\\
Y^R_{2,j}\ &&=&\ \tilde{y}^R_{2,0}\ +\ \sum_{i=1}^{j}\left(1-X_{n-t_n+i}\right)D_{2,n-t_n+i}\\
&&=&\ (1-\tilde{z}^-_{0})\tilde{y}_0\ +\ \sum_{i=n-t_n+1}^{n-t_n+j}\left(1-X_{i}\right)D_{2,i}.
\end{aligned}
\end{equation}
Hence, combining~\eqref{eq:urn_tilde_coupled} and~\eqref{eq:urn_RRU_coupled} with $j=t_n$,
we have that on the set $\tilde{T}_n$
$$(\tilde{Y}^-_{1,n},\tilde{Y}^-_{2,n})=(Y^{R}_{1,t_n},Y^{R}_{2,t_n}).$$
Now, from the asymptotic behavior of the RRU studied in~\cite{Muliere.et.al.06} we have that (since $m_1>m_2$) $\bm{P}(\lim_{n\rightarrow\infty} Z^{R}_n=1)=1$.
Thus, on the set $\tilde{T}_n$ we have $\{\lim_{n\rightarrow\infty} Z^R_n=1\}$, which implies
$\bm{P}\left(H^-_n\cup\tilde{T}_n\right)\rightarrow0$.
This concludes the proof.
\endproof

\proof[Theorem~\ref{thm:exponential_bounds_similar_MRRU_su}]
First, consider the set $F_n=\{y_0+c_1n\leq Y_n \leq y_0+C_1n\}$ defined in~\eqref{def:F_n} and
by using~\eqref{eq:Y_exp_decay} in Theorem~\ref{thm:Y_bounded_and_diverge} we have
$$\overline{\lim}_n\bm{P}(F_n^c)\ =\ 0.$$
Hence, since $|G(n,s_{\delta})|\leq\max\{Z_{n+ns_{\delta}};Z_{n}\}<1$ a.s.,
to prove~\eqref{eq:exponential_bounds_similar_MRRU_su} it is enough to show that for any $0<h<1/2$
\begin{equation}\label{eq:limsup_first_step}
\bm{E}\left[G_{n,s_{\delta}}\ind_{Q^c(\delta,n)\cap F_n}\right]\ \leq\ h\ +\ o(1),
\end{equation}
where we recall that $G(n,s_{\delta})=(\Delta_{n+ns_{\delta}}-\Delta_n)$ and $Q(\delta,n)=\{\Delta_n>\delta\}$.
Now, define $H:=[\delta/h]$ and note that
$$[0,\delta]\ \subset\ [0,(H+1)h]\ =\ \cup_{i=0}^{H}[ih,(i+1)h];$$
then, calling
\begin{equation}\label{def:delta_Q}
\bar{Q}((i+1)h,n):=Q^c((i+1)h,n)\setminus Q^c(ih,n)=\{ih<\Delta_n<(i+1)h\},
\end{equation}
(where for any two sets $A$ and $B$, $A\setminus B=A\cap B^c$),
we have $Q^c(\delta,n)=\cup_{i=0}^{H} \bar{Q}((i+1)h,n)$ and hence the left-hand side of~\eqref{eq:limsup_first_step} can be written as
$$\bm{E}\left[G(n,s_{\delta})\ind_{Q^c(\delta,n)\cap F_n}\right]=
\sum_{i=0}^{H}\bm{E}\left[G(n,s_{\delta})\ind_{\bar{Q}((i+1)h,n)\cap F_n}\right];$$
thus, result~\eqref{eq:limsup_first_step} can be achieved by establishing the following
\begin{equation}\label{eq:limsup_second_step}
\bm{E}\left[G(n,s_{\delta})\ind_{\bar{Q}((i+1)h,n)\cap F_n}\right]\
\leq\ h\cdot\bm{P}\left(\bar{Q}((i+1)h,n)\right)\ +\ o(1),\end{equation}
for any $i\in\{1,..,H\}$.
Now, fix $i\in\{0,..,H\}$, call $\nu:=(i+1)h$ and consider the set $\mathcal{M}_{n}^{\epsilon}:=\{M_n<ns_{\delta}\epsilon\}$ defined in~\eqref{def:M_n_epsilon},
where we recall that $M_n=\sum_{i=n}^{n+ns_{\delta}}\ind_{R^c(\nu,n)}$.
The left-hand side of~\eqref{eq:limsup_second_step} can be so decompose
$\bm{E}\left[G(n,s_{\delta})\ind_{\bar{Q}(\nu,n)\cap F_n}\right]=\mathcal{G}_{1n}+\mathcal{G}_{2n}$, where
$$\mathcal{G}_{1n}:=\bm{E}\left[G(n,s_{\delta})\ind_{\bar{Q}(\nu,n) \cap F_n\cap\mathcal{M}_n^{\epsilon}}\right],\ \mbox{and}\
\mathcal{G}_{2n}:=\bm{E}\left[G(n,s_{\delta})\ind_{\bar{Q}(\nu,n) \cap F_n\cap\mathcal{M}_n^{\epsilon c}}\right].$$
Since $\bm{P}(R(\nu,n))\rightarrow1$ from~\eqref{eq:rho_convergence_probability}, and by using Markov's inequality we have that
$$\bm{P}(\mathcal{M}_{n}^{\epsilon c})\ \leq\ \epsilon^{-1}\frac{1}{ns_{\delta}}
\sum_{i=n}^{n+ns_{\delta}}\bm{P}\left(R^c(\nu,n)\right)\ \rightarrow\ 0;$$
thus, since $|G(n,s_{\delta})|\leq\max\{Z_{n+ns_{\delta}};Z_{n}\}<1$ a.s., we have $\mathcal{G}_{2n}\rightarrow0$ and hence
result~\eqref{eq:limsup_second_step} can be achieved by establishing the following
\begin{equation}\label{eq:limsup_third_step}
\mathcal{G}_{1n}\ =\ \bm{E}\left[G(n,s_{\delta})\ind_{\bar{Q}(\nu,n) \cap F_n\cap\mathcal{M}_n^{\epsilon}}\right]\
\leq\ h\cdot\bm{P}\left(\bar{Q}(\nu,n)\right)\ +\ o(1),
\end{equation}
where we recall that $\bar{Q}(\nu,n)=\{\nu-h<\Delta_n<\nu\}$.

Now, following the same arguments used to determine $s_{\delta}$ in~\eqref{def:s_delta},
we can fix a value $s_h$ such that
$$\bm{P}\left(\ \left\{|G(n,s_h)|> h/2\right\}\ \cap\ F_n\ \right)\ =\ 0,$$
where we recall that $G(n,s_h)=(\Delta_{n+ns_h}-\Delta_n)$.
Analogously to~\eqref{def:s_delta}, the range of values admissible for $s_h$ is
\begin{equation}\label{def:s_h}
s_{h}\ \in\ \left(\ 0\ ,\ \exp\left(\frac{c_1}{2b}h\right)-1\ \right),
\end{equation}
where we recall that $c_1>0$ is a constant introduce in~\eqref{def:F_n} to define $F_n$.

Now, consider the random time $t_j$ defined in~\eqref{def:t_n_lemma} as the smallest time $k$
such that $Q^c(\nu,n-k)$ occurs, i.e. $n-t_n$ indicates the last time up to $n$ the urn proportion is in the interval $(\rho_1-\nu,\rho_1+\nu)$.
Then, call $\tau_n:=t_{n+ns_{\delta}}$ and note that,
since $bar{Q}(\nu,n)\subset Q^c(\nu,n)$ by definition of $\bar{Q}(\nu,n)$, we have that
$$\bm{P}\left(\tau_{n}\leq ns_{\delta}\ |\ \bar{Q}(\nu,n)\right)\ =\ 1.$$
Hence, define $S_H:=[s_{\delta}/s_h]$ and, assuming wlog that $s_{\delta}=S_Hs_hh+1$, on the set $\bar{Q}(\nu,n)$, consider the partition
$\{0,..,ns_{\delta}\}=\cup_{k=0}^{S_H}\mathcal{T}_{k}^n$,
where $\mathcal{T}_{k}^n:=\{nks_h,..,n(k+1)s_h\}$;
thus, the left-hand side of~\eqref{eq:limsup_third_step} can be decompose as $\mathcal{G}_{1n}=\sum_{k=0}^{S_H}T_{k}^n$,
where for any $k\in\{0,..,S_H\}$
\begin{equation}\label{def:T_k_n}
T_{k}^n\ :=\ \bm{E}\left[G(n,s_{\delta})\ind_{\bar{Q}(\nu,n) \cap
F_n\cap\mathcal{M}_n^{\epsilon}\cap\{\tau_n\in\mathcal{T}_{k}^n\}}\right].\end{equation}
Hence, equation~\eqref{eq:limsup_third_step} can be achieved by establishing the following
\begin{equation}\label{eq:limsup_fourth_step}
T_{k}^n\ \leq\ h\cdot\bm{P}\left(\bar{Q}(\nu,n)\cap\{\tau_n\in\mathcal{T}_{k}^n\}\right)\ +\ o(1),\ \ \ \forall k\in\{0,..,S_H\}.
\end{equation}

First, consider $k=0$ in~\eqref{eq:limsup_fourth_step}.
From the definition of $\tau_n$, we have
\begin{equation}\label{eq:relation0_induction2_in_theorem}
\{\tau_n\in\mathcal{T}_{0}^n\}\subset Q^c(\nu+h,n+ns_{\delta}),\end{equation}
where we recall that $Q^c(\nu+h,n+ns_{\delta})=\{\Delta_{n+ns_{\delta}}<\nu+h\}.$
Hence, using~\eqref{eq:relation0_induction2_in_theorem} in~\eqref{def:T_k_n}, it is immediate to obtain~\eqref{eq:limsup_fourth_step}.

For $k\in\{1,..,S_H\}$ in~\eqref{eq:limsup_fourth_step},
from the definition of $\tau_n$ and $E_{n,k}$ in~\eqref{def:E_nk}, we have that
\begin{equation}\label{eq:relation_induction2_in_theorem}
\{\tau_n\in\mathcal{T}_{k}^n\}\subset E(n,n+n(s_{\delta}-ks_h)),\end{equation}
where we recall $E(n,k)=\cup_{j=n}^{k}Q^c(\nu,j)$.
Hence, we can use Lemma~\ref{lem:induction2} with $l_n=n+n(s_{\delta}-ks_h)$,
to obtain, on the set $\mathcal{M}_{n}^{\epsilon}\cap F_n$, for any
$j\in\{n+n(s_{\delta}-ks_h)+1,..,n+ns_{\delta}\}$
\begin{equation}\label{eq:th_induction2_in_theorem}
Q^c(\nu,n+n(s_{\delta}-ks_h)) \subset\ Q^c(\tilde{\Delta}^{*}_j,j)\ \ \ a.s.,
\end{equation}
where we recall that $Q^c(\nu,j)=\{\Delta_j<\nu\}$ and $Q^c(\tilde{\Delta}^{*}_j,j)=\{\Delta_j<\tilde{\Delta}^{*}_j\}$,
$\tilde{\Delta}^{*}_j=\max\left\{\tilde{\Delta}^{-}_{j},\tilde{\Delta}^{+}_{j}\right\}$,
$\tilde{\Delta}^{-}_{j}=|\rho_1-\tilde{Z}^{-}_{j}|$ and $\tilde{\Delta}^{+}_{j}=|\rho_1-\tilde{Z}^{+}_{j}|$.
In particular, by using~\eqref{eq:th_induction2_in_theorem} and since
$\bar{Q}(\nu,n)\subset Q(\nu-h,n)= \{\Delta_n>\nu-h\}$, from~\eqref{def:T_k_n} we obtain
\begin{equation}\label{eq:T_k_n_intermediate_step}
T_{k}^n\ \leq\ \bm{E}\left[(\tilde{\Delta}^{*}_{n+ns_{\delta}}-\nu+h)\ind_{\bar{Q}(\nu,n)\cap F_n\cap\mathcal{M}_n^{\epsilon}
\cap\{\tau_n\in\mathcal{T}_{k}^n\}}\right].\end{equation}

Note that, from the definition of $\tau_n$ and $\mathcal{T}_{k}^n$, we have
$$\{\tau_n\in\mathcal{T}_{k}^n\}\ \subset\ \{nks_h<t_{n+ns_{\delta}}<n(k+1)s_h\}.$$
Hence, we can apply Lemma~\ref{lem:comparison_tilde_RRU} with $k_{n+ns_{\delta}}=nks_h$,
$\tilde{T}_j:=\left\{\tilde{\Delta}^{*}_j>\nu\right\}$ and
$H_j:=\{k_{j}<t_{j}<\infty\}$ as defined in~\eqref{def:Ttilde_n_H_n}, so obtaining
$$\bm{E}\left[(\tilde{\Delta}^{*}_{n+ns_{\delta}}-\nu)^+\ind_{\{\tau_n\in\mathcal{T}_{k}^n\}}\right]\ \leq\ \bm{P}\left(H_{n+ns_{\delta}}\cup\tilde{T}_{n+ns_{\delta}}\right)\rightarrow0.$$
Hence, applying these results to~\eqref{eq:T_k_n_intermediate_step}, we obtain
$$T_{k}^n\ \leq\ h\cdot \bm{P}\left(\bar{Q}(\nu,n)\cap\{\tau_n\in\mathcal{T}_{k}^n\}\right)\ +\ o(1),$$
that corresponds to~\eqref{eq:limsup_fourth_step}.
This concludes the proof.
\endproof

\subsection{Proof of weak consistency}   \label{subsection_proof_WLLN}

\proof[Theorem~\ref{thm:convergence_probability}]
The result is established by proving that, for any $l>0$ and any $\epsilon>0$, there exists $n_0\in\mathbb{N}$ such that
\begin{equation}\label{eq:th_convergence_WLLN}
\bm{P}\left(Q(l,n)\right)<\epsilon,
\end{equation}
for any $n\geq n_0$, where we recall that $Q(l,n)=\{\Delta_n>l\}$ and $\Delta_n=|\rho_1-Z_n|$.
To this end, fix $0<\epsilon^{\prime}<\frac{l\epsilon}{3}$ and $0<\delta<\epsilon^{\prime}$ to define the conditions
$$\mathcal{A}_n\ :=\ \{\bm{P}\left(Q(\delta,n)\right)<\epsilon^{\prime}\},\ \qquad\
\mathcal{B}_n\ :=\ \{\bm{E}[\Delta_n]<2\epsilon^{\prime}\}.$$
It is immediate to see that $\mathcal{B}_n$ implies~\eqref{eq:th_convergence_WLLN}.
Thus,~\eqref{eq:th_convergence_WLLN} can be established by proving that
\begin{itemize}
\item[(a)] for any $N\geq1$ there exists $n_0\geq N$ such that $\mathcal{A}_{n_0}$ occurs;
\item[(b)] there exists $n_0\geq1$ such that for any $n\geq n_0$ $\mathcal{A}_{n}\subset \mathcal{B}_{k}$ for all $k\in\{n+1,..,n(1+s_{\delta})\}$;
\item[(c)] there exists $n_0\geq1$ such that for any $n\geq n_0$ $\mathcal{B}_{n}\subset \mathcal{B}_{k}$ for all $k\in\{n(1+s_{\delta}),..,(n+1)(1+s_{\delta})\}$.
\end{itemize}
For part (a), we will show that cannot exist $N\geq1$ such that
\begin{equation}\label{eq:contradiction_hp_WLLN}
\mathcal{A}^c_n\ :=\ \{\bm{P}\left(Q(\delta,n)\right)\geq\epsilon^{\prime}\},
\end{equation}
occurs for all $n\geq N$. First, we combine
Theorem~\ref{thm:exponential_bounds_similar_MRRU_giu} and Theorem~\ref{thm:exponential_bounds_similar_MRRU_su} to obtain
\begin{equation}\label{eq:relation_final_convergence}
\bm{E}\left[G(n,s_{\delta})\right]\ \leq\ -C\left(\bm{P}(Q(\delta,n))\ -\ \frac{\epsilon^{\prime}}{2}\right),
\end{equation}
with $0<C<\infty$, where we recall that $G(n,s_{\delta})=(\Delta_{n+ns_{\delta}}-\Delta_n)$.
Now, if~\eqref{eq:contradiction_hp_WLLN} holds, then there exists a subsequence $\{k_n;n\geq1\}$ such that,
$k_1=N$ and $k_n=k_{n-1}(1+s_{\delta})$ for all $n\geq2$, and by~\eqref{eq:relation_final_convergence}
$$\bm{E}\left[\Delta_{k_n}\right]=\sum_{i=1}^n\bm{E}\left[G(k_{i-1},s_{\delta})\right]
\leq -\sum_{i=1}^n C\frac{\epsilon^{\prime}}{2}\ =\ -\infty,$$
where $G(k_{i-1},s_{\delta})=(\Delta_{k_i}-\Delta_{k_{i-1}})$,
which is a contradiction and hence part (a) holds.
For part (b), consider the time $n$ at which $\mathcal{A}_n$ occurs.
Fix $k\in\{n+1,..,n+ns_{\delta}\}$ and note that $\bm{E}[\Delta_k]\ \leq\ J_{1n}\ +\ J_{2n,k}$ where
$$J_{1n}\ :=\ \bm{E}[\Delta_n],\ \ \mbox{and}\ \ J_{2n,k}\ :=\ \bm{E}[|\Delta_k-\Delta_n|].$$
From definition of $s_{\delta}$ in~\eqref{def:s_delta} we have
\[\begin{aligned}
J_{2n,k}\ &&\leq&\ \bm{E}[|\Delta_k-\Delta_n|\ind_{F_n}]\ +\ \bm{E}[|\Delta_k-\Delta_n|\ind_{F^c_n}]\\
&&\leq&\ \delta\ +\ \bm{P}(F^c_n),
\end{aligned}\]
and using $\bm{P}(F_n^c)\rightarrow0$ from~\eqref{eq:Y_exp_decay} in Theorem~\ref{thm:Y_bounded_and_diverge} we have that
$\lim_{n\rightarrow\infty}J_{2n,k}\leq \delta$. Thus,
there exists $n_0\geq1$ such that $J_{2n,k}<2\delta$ for any $n\geq n_0$.
Then, note that $J_{1n}\ =\ J_{3n}\ +\ J_{4n}$ where
$$J_{3n}\ :=\ \bm{E}[\Delta_n\ind_{Q^c(\delta,n)}],\ \ \mbox{and}\ \ J_{4n}\ :=\ \bm{E}[\Delta_n\ind_{Q(\delta,n)}].$$
Notice that $J_{3n}\leq\delta\bm{P}(Q^c(\delta,n))<\delta$ and
$J_{4n}\leq\bm{P}(Q(\delta,n))<\epsilon^{\prime}$, and hence we have $J_{1n}<\delta+\epsilon^{\prime}$.
Thus, combining $J_{1n}$ and $J_{2n}$, since $\delta<\epsilon^{\prime}/3$, we obtain for any $n\geq n_0$
$$\bm{E}[\Delta_k]\ \leq\ J_{1n}\ +\ J_{2n,k}\ < \delta\ +\ \epsilon^{\prime}\ +\ 2\delta\ <\ 2\epsilon^{\prime},$$
that implies (b).
For part (c), for any $k\in\{n(1+s_{\delta}),..,(n+1)(1+s_{\delta})\}$ consider
$$\bm{E}[|\Delta_k-\Delta_{n+ns_{\delta}}|]\ \leq\ \bm{E}[|\Delta_k-\Delta_{n+ns_{\delta}}|\ind_{F_n}]\ +\ \bm{E}[|\Delta_k-\Delta_{n+ns_{\delta}}|\ind_{F^c_n}].$$
First, note that $\bm{P}(F_n^c)\rightarrow0$ from~\eqref{eq:Y_exp_decay} in Theorem~\ref{thm:Y_bounded_and_diverge}.
Then, since $|k-(n+ns_{\delta})|\leq(1+s_{\delta})$ and $|Z_{n+1}-Z_n|<b/Y_n$ a.s., we have that
$$\bm{P}\left(\ \left\{|Z_k-Z_{n+ns_{\delta}}|>\left(\frac{b}{y_0+c_1n}\right)(1+s_{\delta})\right\}\ \cap\ F_n\ \right)\ =\ 0.$$
Thus, for any $k\in\{n(1+s_{\delta}),..,(n+1)(1+s_{\delta})\}$ we have
\begin{equation}\label{eq:delta_k_delta_n}
\bm{E}[|\Delta_k-\Delta_{n+ns_{\delta}}|]\ \leq\ \left(\frac{b(1+s_{\delta})}{y_0+c_1n}\right)\ +\ \bm{P}(F_n^c)\ \rightarrow\ 0.\end{equation}
Now, since $\mathcal{B}_{n}\subset\mathcal{A}_{n}\cup\mathcal{C}_{n}$, where $\mathcal{C}_{n}=(\mathcal{B}_{n}\cap\mathcal{A}^c_{n})$,
part (c) is established by proving that there exists $n_0\geq1$ such that, for any $n\geq n_0$,
\begin{itemize}
\item[(c1)] $\mathcal{A}_{n}\subset \mathcal{B}_{k}$ for all $k\in\{n(1+s_{\delta}),..,(n+1)(1+s_{\delta})\}$;
\item[(c2)] $\mathcal{C}_{n}\subset \mathcal{B}_{k}$ for all $k\in\{n(1+s_{\delta}),..,(n+1)(1+s_{\delta})\}$;
\end{itemize}
For part (c1), we can follow the same arguments of part (b), except for $J_{2n,k}$ since here $k\in\{n(1+s_{\delta}),..,(n+1)(1+s_{\delta})\}$ and hence
\[\begin{aligned}
J_{2n,k}\ &&\leq&\ \bm{E}[|\Delta_k-\Delta_{n}|\ind_{F_n}]\ +\ \bm{E}[|\Delta_k-\Delta_n|\ind_{F^c_n}]\\
&&\leq&\ \bm{E}[|\Delta_k-\Delta_{n+ns_{\delta}}|\ind_{F_n}]\ +\ \bm{E}[|\Delta_{n+ns_{\delta}}-\Delta_n|\ind_{F_n}]\ +\ \bm{P}(F^c_n)\\
&&\leq&\ \bm{E}[|\Delta_k-\Delta_{n+ns_{\delta}}|\ind_{F_n}]\ +\ \delta\ +\ \bm{P}(F^c_n);
\end{aligned}\]
However, by using~\eqref{eq:delta_k_delta_n}, we still have $\lim_{n\rightarrow\infty}J_{2n,k}\leq \delta$ and so, analogously to part (b),
there exists $n_0\geq1$ such that $J_{n2}<2\delta$ for any $n\geq n_0$.
Since $J_{1n}$ does not depend on $k$, (c1) follows.
For part (c2), we combine~\eqref{eq:relation_final_convergence} and $\mathcal{A}^c_{n}$ to obtain
\begin{equation}\label{eq:relation_final_convergence_2}
\bm{E}\left[G(n,s_{\delta})\right]\ \leq\ -C\frac{\epsilon^{\prime}}{2}\,
\end{equation}
where we recall that $G(n,s_{\delta})=(\Delta_{n+ns_{\delta}}-\Delta_n)$.
Moreover, by~\eqref{eq:delta_k_delta_n} there exists $n_0\geq1$ such that $\bm{E}[|\Delta_k-\Delta_{n+ns_{\delta}}|]\leq C\frac{\epsilon^{\prime}}{2}$
for any $n\geq n_0$.
Hence, (c2) follows by combining~\eqref{eq:delta_k_delta_n},~\eqref{eq:relation_final_convergence_2} and $\mathcal{B}_n$ as follows:
$$\bm{E}[\Delta_k]\ \leq\ \bm{E}[|\Delta_k-\Delta_{n+ns_{\delta}}|]\ +\ \bm{E}[G(n,s_{\delta})]\ +\ \bm{E}[\Delta_n]\ =\ 2\epsilon^{\prime}.$$
\endproof

\begin{rem}
It is possible to present a modification of the current arguments along the traditional probabilistic lines.
We chose to present the above alternative logical argument.
\end{rem}

\section{Proof of strong consistency}   \label{section_SLLN}

In this section, we provide the proof of the strong consistency of the urn proportion $Z_n$ for any values of $m_1$ and $m_2$,
when the random thresholds $\hat{\rho}_{1,n}$ and $\hat{\rho}_{2,n}$ converge with probability one.

\proof[Theorem~\ref{thm:convergence_as}]
We divide the proof in three steps:
\begin{description}
\item[(a)] $\bm{P}\left(\ \rho_2\ \leq\ \underline{\lim}_nZ_n \leq \overline{\lim}_nZ_n\ \leq\ \rho_1\ \right)\ =\ 1$,
\item[(b)]
\begin{equation*}
\begin{cases}
\bm{P}(\overline{\lim}_nZ_n\geq\rho_1)=1\ &\text{if } m_1> m_2,\\
\bm{P}(\underline{\lim}_nZ_n\leq\rho_2)=1\ &\text{if } m_1< m_2.
\end{cases}
\end{equation*}
\item[(c)] $\bm{P}\left(\ \lim_nZ_n \textit{ exists}\ \right)\ =\ 1.$
\end{description}

For part (a), firstly note that, when $\rho_1=1$ and $\rho_2=0$, result (a) is trivially true, hence consider $0<\rho_2\leq\rho_1<1$.
We show that $\bm{P}(\overline{\lim}_nZ_n \leq\rho_1)=1$, since the proof of
$\bm{P}(\underline{\lim}_nZ_n \geq\rho_2)=1$ is completely analogous.
To this end, we show that cannot exist $\epsilon>0$ and $\rho^{\prime} >\rho_1$ such that
\begin{equation}\label{eq:hp_assurda_part_a}
\bm{P}\left(\overline{\lim}_nZ_n>\rho_1^{\prime}\right)\ \geq\ \epsilon\ >\ 0.
\end{equation}
We prove this by contradiction using a comparison argument with an RRU model.
The proof involves last exit time arguments.
Now, suppose~\eqref{eq:hp_assurda_part_a} holds and let $A_1:=\{\overline{\lim}_nZ_n>\rho_1^{\prime}\}$.
Let
$$R_1:=\left\{\ k\geq0\ :\ \hat{\rho}_{1,k}\geq \frac{\rho_1^{\prime}+\rho_1}{2}\ \right\},$$
and denote the last time the process $\{\hat{\rho}_{1,n};n\geq1\}$ is above $\left(\rho_1^{\prime}+\rho_1\right)/2$ by
\[\begin{aligned}
t_{\frac{\rho_1^{\prime}+\rho_1}{2}} & =
\begin{cases}
\sup\{R_1\} & \text{if }R_1 \neq\emptyset;\\
0 & \text{otherwise}.
\end{cases}
\end{aligned}\]
Since $\hat{\rho}_{1,n}\stackrel{a.s.}{\rightarrow}\rho_1$ by~\eqref{eq:rho_convergence_as}, then we have that $\bm{P}\left(t_{\frac{\rho_1^{\prime}+\rho_1}{2}}<\infty\right)=1$.
Hence, there exists $n_{\epsilon}\in\mathbb{N}$ such that
\begin{equation}\label{eq:epsilon_mezzi_1}
\bm{P}\left(t_{\frac{\rho_1^{\prime}+\rho_1}{2}}>n_{\epsilon}\right)\ \leq\ \frac{\epsilon}{2}.
\end{equation}
Setting $B_1:=\left\{t_{\frac{\rho_1^{\prime}+\rho_1}{2}}>n_{\epsilon}\right\}$ and using~\eqref{eq:epsilon_mezzi_1}, it follows that
$$\epsilon\ \leq\ \bm{P}\left(A_1\right)\ \leq\ \epsilon/2\ +\ \bm{P}\left(A_1\cap B_1^c\right).$$
Now, we show that $\bm{P}\left(A_1\cap B_1^c\right)=0$.
Setting
$$C_1\ =\ \left\{\ \omega\in\Omega\ :\ \underline{\lim}_nZ_n<\frac{\rho_1^{\prime}+\rho_1}{2}\ \right\},$$
we decompose $\bm{P}\left(A_1\cap B_1^c\right)$ as follows:
$$\bm{P}\left(A_1\cap B_1^c\right)\ \leq\ \bm{P}\left(E_1\right)\ +\ \bm{P}\left(E_2\right),$$
where $E_1=A_1\cap B_1^c\cap C_1$ and $E_2=A_1\cap B_1^c\cap C_1^c.$

Consider the term $\bm{P}\left(E_2\right)$. Note that on the set $C_1^c$, we have $\left\{\underline{\lim}_nZ_n\geq\frac{\rho_1^{\prime}+\rho_1}{2} \right\}$ and  on the set $B_1^c$ we have $\{\hat{\rho}_{1,n}\leq \frac{\rho_1^{\prime}+\rho_1}{2}\}$ for any $n\geq n_{\epsilon}$.
Hence, since $B_1^c\cap C_1^c\supset E_2$, on the set $E_2$ we have that $W_{1,n}=\ind_{\{Z_n\leq \hat{\rho}_{1,n}\}}\stackrel{a.s.}{\rightarrow} 0$.
Then, letting $\tau_W:=\sup\{k\geq1:W_{1,k}=1\}$ we have $\bm{P}(E_2\cap\{\tau_W<\infty\})=\bm{P}(E_2)$ and,
on the set $E_2$, for any $n\geq \tau_W$ the ARRU model can be written as follows:
\[\left\{
\begin{array}{l}
Y_{1,n+1}=Y_{1,\tau_W}\\
\\
Y_{2,n+1}=Y_{2,\tau_W}+\sum_{i=\tau_W}^{n+1}\left(1-X_i\right)D_{2,i},
\end{array}
\right.\]
where $W_{1,i-1}=0$ for any $i\geq\tau_W$, and $W_{2,i-1}=1$ because $W_{2,i-1}+W_{2,i-1}\geq1$ by construction.
Now, consider an RRU model $\{Z^{R}_{i};i\geq1\}$ with initial composition $(Y^{R}_{1,0},Y^{R}_{2,0})=(Y_{1,\tau_W},Y_{2,\tau_W})$ a.s.;
the reinforcements are defined as $D^{R}_{1,i}=0$ and $D^{R}_{2,i}=D_{2,\tau_W+i}$ for any $i\geq1$ a.s.;
the drawing process is modeled by $X^R_{i+1}:=\ind_{\{U^R_{i+1}<Z^R_{i}\}}$ and $U^{R}_{i}=U_{\tau_W+i}$ a.s.,
where $\{U_n;n\geq1\}$ is the sequence such that $X_{n+1}=\ind_{\{U_n<Z_n\}}$ for any $n\geq1$.
Formally, this RRU model can be described for any $n\geq1$ as follows:
\[\left\{
\begin{array}{l}
Y^R_{1,n+1}=Y^R_{1,0}=Y_{1,\tau_W}\\
\\
Y^R_{2,n+1}=Y^R_{2,0}+\sum_{i=0}^{n+1}\left(1-X^R_i\right)D^R_{2,i}=Y_{2,\tau_W}+\sum_{i=\tau_W}^{n+\tau_W+1}\left(1-X_{i}\right)D_{2,i}.
\end{array}
\right.\]
Hence, on the set $E_2$ we have that for any $n\geq\tau_W$
$$(Y_{1,n},Y_{2,n})=(Y^{R}_{1,n-\tau_W},Y^{R}_{2,n-\tau_W}).$$
Since from~\cite{Muliere.et.al.06} $\bm{P}(\overline{\lim}_n Z^{R}_n=0)=1$,
on the set $E_2$ we have that $\{\overline{\lim}_n Z_n=0\}$.
This is incompatible with the set $A_1$ which includes $E_2$. Hence $\bm{P}\left(E_2\right)=0$.

We now turn to the proof that $\bm{P}\left(E_1\right)=0$.
To this end, let
$$\tau_{\epsilon}\ :=\ \inf\left\{k\geq n_{\epsilon}\ :\ \left\{Z_k<\frac{\rho_1^{\prime}+\rho_1}{2}\right\}\cap
\left\{Y_k>\frac{b}{(\rho_1^{\prime}-\rho_1)/2}\right\}\ \right\}$$
and note that, since by Lemma~\ref{lem:urn_not_stop} $Y_n\stackrel{a.s.}{\rightarrow}\infty$,
$\bm{P}(C_1\cap\{\tau_{\epsilon}<\infty\})=\bm{P}(C_1)$.
Moreover, on the set $B_1^c$ we have that $\{\hat{\rho}_{1,n}\leq \frac{\rho_1^{\prime}+\rho_1}{2}\}$ for any $n\geq n_{\epsilon}$.
We now show by induction that on the set $B_1^c\cap C_1$ we have $\{Z_{n}<\rho_1^{\prime}\ \forall n\geq \tau_{\epsilon}\}$.
By definition we have $Z_{\tau_{\epsilon}}<\frac{\rho_1^{\prime}+\rho_1}{2}$, and by Lemma~\ref{lem:Y_increments} this implies $Z_{\tau_{\epsilon}+1}<\rho_1^{\prime}$;
now, consider an arbitrary $n> \tau_{\epsilon}$; if $Z_n<\frac{\rho_1^{\prime}+\rho_1}{2}$, then by Lemma~\ref{lem:Y_increments} we have $Z_{n+1}<\rho_1^{\prime}$;
if $\frac{\rho_1^{\prime}+\rho_1}{2}<Z_n<\rho_1^{\prime}$ we have $W_{1,n}=0$ and so $Z_{n+1}\leq Z_n<\rho_1^{\prime}$.
Hence, since $B_1^c\cap C_1\subset E_1$, on the set $E_1$ we have $\{Z_{n}<\rho_1^{\prime}\ \forall n\geq \tau_{\epsilon}\}$.
This is incompatible with the set $A_1$ which also includes $E_1$. Hence $\bm{P}\left(E_1\right)=0$.
Combining all together we have $\epsilon\ \leq\ \epsilon/2 + \bm{P}\left(E_1\right) + \bm{P}\left(E_2\right)\ =\ \epsilon/2$,
which is impossible.
Thus, we conclude that $\bm{P}(A_1^c)=\bm{P}(\overline{\lim}_nZ_n \leq\rho_1)=1$.

For part (b), wlog we assume $m_1>m_2$ to show that $\bm{P}(\overline{\lim}_nZ_n \geq\rho_1)=1$,
since the proof of $\bm{P}(\underline{\lim}_nZ_n \leq\rho_2)=1$ when $m_1<m_2$ is completely analogous.
To this end, we now show that cannot exist $\epsilon>0$ and $\rho^{\prime} <\rho_1$ such that
\begin{equation}\label{eq:hp_assurda_part_b}
\bm{P}\left(\overline{\lim}_nZ_n<\rho_1^{\prime}\right)\ \geq\ \epsilon\ >\ 0.
\end{equation}
We prove this by contradiction, using a comparison argument with an RRU model.
Now suppose~\eqref{eq:hp_assurda_part_b} holds and let $A_2:=\{\overline{\lim}_nZ_n<\rho_1^{\prime}\}$.
Let
$$R_2:=\left\{\ k\geq0\ :\ \hat{\rho}_{1,k}< \frac{\rho_1^{\prime}+\rho_1}{2}\ \right\},$$
and define the last time the process $\{\hat{\rho}_{1,n};n\geq1\}$ is less than $\left(\rho_1^{\prime}+\rho_1\right)/2$ by
\[\begin{aligned}
\tau_{\frac{\rho_1^{\prime}+\rho_1}{2}} & =
\begin{cases}
\sup\{R_2\} & \text{if }R_2 \neq\emptyset;\\
0 & \text{otherwise}.
\end{cases}
\end{aligned}\]
Since $\hat{\rho}_{1,n}\stackrel{a.s.}{\rightarrow}\rho_1$, then we have that $\bm{P}\left(\tau_{\frac{\rho_1^{\prime}+\rho_1}{2}}<\infty\right)=1$.
Hence, there exists $n_{\epsilon}\in\mathbb{N}$ such that
\begin{equation}\label{eq:epsilon_mezzi_2}
\bm{P}\left(\tau_{\frac{\rho_1^{\prime}+\rho_1}{2}}>n_{\epsilon}\right)\ \leq\ \frac{\epsilon}{2}.
\end{equation}
Setting $B_2:=\left\{\tau_{\frac{\rho_1^{\prime}+\rho_1}{2}}>n_{\epsilon}\right\}$ and using~\eqref{eq:epsilon_mezzi_2}, it follows that
$$\epsilon\ \leq\ \bm{P}\left(A_2\right)\ \leq\ \epsilon/2\ +\ \bm{P}\left(A_2\cap B_2^c\right).$$
Let $E_3:=A_2\cap B_2^c$. We now show that $\bm{P}\left(E_3\right)$=0.
On the set $A_2$, we have $\left\{\underline{\lim}_n Z_n\leq\rho_1^{\prime} \right\}$ and
on the set $B_2^c$,
we have $\{\hat{\rho}_{1,n}\geq \frac{\rho_1^{\prime}+\rho_1}{2}\}$ for any $n\geq n_{\epsilon}$.
Hence, on the set $E_3$ we have that $W_{1,n}=\ind_{\{Z_n\leq \hat{\rho}_{1,n}\}}\stackrel{a.s.}{\rightarrow} 1$.
Then, letting $\tau_W:=\sup\{k\geq1:W_{1,n}=0\}$ we have $\bm{P}(E_3\cap\{\tau_W<\infty\})=\bm{P}(E_3)$.
Now, analogously to the proof of $\bm{P}\left(E_2\right)=0$, we can use comparison arguments with the RRU model
to show that on the set $E_3$ we have $\{\overline{\lim}_nZ_n=1\}$.
This is incompatible with the set $A_2$, which also includes $E_3$. Hence $\bm{P}\left(E_3\right)=0$.
Combining all together we have $\epsilon\ \leq\ \epsilon/2 + \bm{P}\left(E_3\right)\ =\ \epsilon/2$,
which is impossible.
Thus, we conclude that the event $A_2^c=\{\overline{\lim}_nZ_n \geq\rho_1\}$ occurs with probability one.

For part (c), note that, combining (a) and (b), we have shown that
\begin{equation}
\begin{cases}
\bm{P}(\overline{\lim}_nZ_n=\rho_1)=1\ &\text{if } m_1> m_2,\\
\bm{P}(\rho_2\leq\underline{\lim}_nZ_n\leq\overline{\lim}_nZ_n\leq\rho_1)=1\ &\text{if } m_1= m_2,\\
\bm{P}(\underline{\lim}_nZ_n=\rho_2)=1\ &\text{if } m_1< m_2.
\end{cases}
\end{equation}
Therefore, if the process $\{Z_n;n\geq1\}$ converges almost surely, we obtain~\eqref{eq:Z_convergence_as}.
Wlog, assume $m_1\geq m_2$, since the proof of the case $m_1\leq m_2$ is completely analogous.

First, let $d$, $u$, $\gamma$ and $\rho_1^{\prime}$ ($d<u<\gamma<\rho_1^{\prime}<\rho_1$) be four constants in $\left(0,1\right)$.
Let $\{\tau_j(d,u);j\geq1\}$ and $\{t_j(d,u);j\geq1\}$ be the sequences of random variables defined in~\eqref{eq:def_tau}.
Since $d$ and $u$ are fixed in this proof, we sometimes denote $\tau_j(d,u)$ by $\tau_j$ and $t_j(d,u)$ by $t_j$.
It is easy to see that $\tau_n$ and $t_n$ are stopping times with respect to $\left\{\mathcal{F}_n;n\geq1\right\}$.

Recall that, by Lemma~\ref{lem:Zn_converge_ifonlyif}, we have that for every $0<d<u<1$
\[\begin{aligned}
Z_n\ \text{converges a.s.}\ &&\Leftrightarrow&\ \bm{P}\left( t_{n}(d,u) < \infty \right)\rightarrow0,\\
&&\Leftrightarrow&\ \sum_{n=1}^{\infty} \bm{P}\left( t_{n+1}(d,u) = \infty | t_n(d,u) < \infty\right) =\infty.
\end{aligned}\]
Now, to prove that $Z_n$ converges a.s., it is sufficient to show that
$$\bm{P}\left( t_n(d,u) < \infty \right)\rightarrow0,$$
for all $0<d<u<1$.
Suppose $Z_n$ does not converges a.s.. This implies that $\bm{P}\left( t_n < \infty \right)\downarrow \phi_1>0$,
since $\bm{P}\left( t_n < \infty \right)$ is a non-increasing sequence.
We will show that for large $j$ there exists a constant $\phi<1$ dependent on $\phi_1$, such that
\begin{equation}\label{eq:goal_a_s_convergence}
\bm{P}\left( t_{j+1} < \infty | t_j<\infty \right) \leq \phi.
\end{equation}
This result implies that $\sum_{n} \bm{P}\left( t_{n+1} = \infty | t_n < \infty\right) =\infty$,
establishing by Lemma~\ref{lem:Zn_converge_ifonlyif} that $\bm{P}\left( t_{n} < \infty \right)$ converges to zero as $n$ goes to infinity,
which is a contradiction.\\

Consider the term $\bm{P}\left( t_{i+1} < \infty | t_i < \infty\right)$.
First, let us denote by $\tau_{\rho_1^{\prime}}$ the last time the process $\hat{\rho}_{1,n}$ is below $\rho_1^{\prime}$, i.e.
\begin{equation*}
\begin{aligned}
\tau_{\rho_1^{\prime}} & =
\begin{cases}
\sup\{n\geq 1 :\hat{\rho}_{1,n}\leq\rho_1^{\prime}\} &
\text{if } \{n\geq 1 :\hat{\rho}_{1,n}\leq\rho_1^{\prime}\}  \neq\emptyset;
\\
0 & \text{otherwise}.
\end{cases}
\end{aligned}
\end{equation*}
Since $\hat{\rho}_{1,n}\stackrel{a.s.}{\rightarrow}\rho_1$, we have that $\bm{P}\left(\tau_{\rho_1^{\prime}}<\infty\right)=1$.
Hence, for any $\epsilon\in\left(0,\frac{1}{2}\right)$ there exists $n_{\epsilon}\in\mathbb{N}$ such that
\begin{equation}\label{eq:epsilon_convergence}
\frac{1}{\phi_1}\bm{P}\left(\tau_{\rho_1^{\prime}}> n_{\epsilon}\right)\ \leq\ \epsilon.
\end{equation}
By denoting $\bm{P_{i}}\left(\cdot\right)=\bm{P}\left(\cdot|t_i<\infty\right)$ and
using $t_i\leq\tau_i\leq t_{i+1}$ we obtain
$$\bm{P}\left(t_{i+1}<\infty|t_i<\infty\right)\ \leq\ \bm{P_{i}}\left(\tau_i<\infty\right).$$
Hence
\begin{equation}\label{eq:decomposition_as}
\bm{P_{i}}\left(\tau_i<\infty\right)\leq\bm{P_{i}}\left(\{\tau_i<\infty\}\cap\{\tau_{\rho_1^{\prime}}\leq n_{\epsilon}\}\right)\ +\
\bm{P_{i}}\left(\tau_{\rho_1^{\prime}}> n_{\epsilon}\right).
\end{equation}
We start with the second term in~\eqref{eq:decomposition_as}. Note that
$$\bm{P_{i}}\left(\tau_{\rho_1^{\prime}}> n_{\epsilon}\right)\ \leq\ \frac{\bm{P}\left(\tau_{\rho_1^{\prime}}> n_{\epsilon}\right)}{\bm{P}\left(t_i<\infty\right)}\ \leq\
\frac{\bm{P}\left(\tau_{\rho_1^{\prime}}> n_{\epsilon}\right)}{\phi_1}\ \leq\ \epsilon,$$
where the last inequality follows from~\eqref{eq:epsilon_convergence}.

Now, consider the first term in~\eqref{eq:decomposition_as}.
Since the probability is conditioned to the set $\{t_i<\infty\}$,
in what follows we will consider the urn process at times $n$ after the stopping time $t_i$.
Since we want to show~\eqref{eq:goal_a_s_convergence} for large $i$,
we can choose an integer $i\geq n_{\epsilon}$ and
$$i> \log_{\frac{u\left(1-d\right)}{d\left(1-u\right)}}\left(\frac{b}{Y_0\left(\gamma-u\right)}\right),$$
so that
\begin{itemize}
\item[(i)] $t_i\geq i\geq n_{\epsilon}$ a.s.;
\item[(ii)] from Lemma~\ref{lem:Y_geometric_increasing}, we have that $Y_{\tau_i} > b/\left(\gamma-u\right)$ a.s.
\end{itemize}
These two properties imply respectively that, on the set $\{n\geq t_i\}$
\begin{itemize}
\item[(i)] $\hat{\rho}_{1,n}\geq\rho_1^{\prime}$, since from $\{\tau_{\rho_1^{\prime}}\leq n_{\epsilon}\}$ we have that $n\geq \tau_{\rho_1^{\prime}}$;
\item[(ii)] $Z_{t_i}\in\left(u,\gamma\right)$, since $Z_{t_i-1}\leq u$ and $Z_{t_i}> u$ and
from Lemma~\ref{lem:Y_increments} we have that $|Z_{n}-Z_{n-1}|<\left(\gamma-u\right)$.
\end{itemize}

Now, let us define two sequences of stopping times $\{t^\ast_n;n\geq1\}$ and $\{\tau^\ast_n;n\geq1\}$,
where $t^\ast_n$ represents the first time after $\tau^\ast_{n-1}$ the process
$Z_{t_i+n}$ up-crosses $\rho_1^{\prime}$,
while $\tau^\ast_n$ represents the first time after $t^\ast_{n}$ the process
$Z_{t_i+n}$ down-crosses $\gamma$.
Formally, let $\tau^\ast_{0}=0$ and define for every $j\geq 1$ two stopping times
\begin{equation}\label{eq:def_tauast2}
\begin{aligned}
t^\ast_j & =
\begin{cases}
\inf\{n> \tau^\ast_{j-1}:Z_{t_i+n}>\rho_1^{\prime}\} &
\text{if } \{n>\tau^\ast_{j}:Z_{t_i+n}>\rho_1^{\prime}\} \neq\emptyset;
\\
+ \infty & \text{otherwise}.
\end{cases}
\\
\tau^\ast_j & =
\begin{cases}
\inf\{n>t^\ast_{j} :Z_{t_i+n}\leq\gamma\} &\ \ \ \
\text{if } \{n>t^\ast_{j-1} :Z_{t_i+n}\leq\gamma\} \neq\emptyset;
\\
+ \infty &\ \ \ \ \text{otherwise}.
\end{cases}
\end{aligned}
\end{equation}
Note that, since $Z_{t_i+\tau^\ast_j-1} \geq \gamma$ and $Z_{t_i+\tau^\ast_j} < \gamma$,
from (ii) we have that $Z_{t_i+\tau^\ast_j}\in\left(u,\gamma\right)$.\\

\noindent For any $j\geq0$, let $\{\tilde{Z}^j_n;n\geq1\}$ be an RRU model defined as follows:
\begin{itemize}
\item[(1)] $\left(\tilde{Y}^j_{1,0},\tilde{Y}^j_{2,0}\right) = \left(Y_{1,t_i+\tau^\ast_j},Y_{1,t_i+\tau^\ast_j}\frac{u+d}{2-u-d}\right)$ a.s.,
which implies that $\tilde{Z}^j_0=\frac{u+d}{2}$;
\item[(2)] the drawing process is modeled by
$\tilde{X}^j_{n+1}=\ind_{\{\tilde{U}^j_{n+1}<\tilde{Z}^j_n\}}$, where $\tilde{U}^j_{n+1}=U_{t_i+\tau^\ast_j+n+1}$ a.s.
and $U_n$ is such that $X_{n}=\ind_{\{U_{n}<Z_{n-1}\}}$;
\item[(3)] the reinforcements are defined as
$\tilde{D}^j_{2,n+1}=D_{2,t_i+\tau^\ast_j+n+1}+\left(m_1-m_2\right)$,
$\tilde{D}^j_{1,n+1}=D_{1,t_i+\tau^\ast_j+n+1}$ a.s.;
this means $\bm{E}[\tilde{D}^j_{1,n}]=\bm{E}[\tilde{D}^j_{2,n}]$ for any $n\geq1$;
\item[(4)] the urn process evolves as an RRU model, i.e. for any $n\geq0$
\[\left\{
\begin{array}{l}
\tilde{Y}^j_{1,n+1}=\tilde{Y}^j_{1,n}+\tilde{X}^j_{n+1}\tilde{D}^j_{1,n+1},\\
\tilde{Y}^j_{2,n+1}=\tilde{Y}^j_{2,n}+\left(1-\tilde{X}^j_{n+1}\right)\tilde{D}^j_{2,n+1},\\
\tilde{Y}^j_{n+1}=\tilde{Y}^j_{1,n+1}+\tilde{Y}^j_{2,n+1},\\
\tilde{Z}^j_{n+1}=\frac{\tilde{Y}^j_{1,n+1}}{\tilde{Y}^j_{n+1}}.
\end{array}
\right.\]
\end{itemize}
We will compare the process $\{\tilde{Z}^j_{n};n\geq1\}$ with the ARRU process $\{Z_{t_i+n};n\geq1\}$.
Note that at time $n$, we have defined only the processes $\tilde{Z}^j$ such that $\tau^\ast_j<n$.\\

\noindent We will prove, by induction, that on the set $\{\tau_{\rho_1^{\prime}}\leq n_{\epsilon}\}$, for any $j\in\mathbb{N}$
and for any $n\leq {t^\ast_{j+1}}-{\tau^\ast_j}$
\begin{equation}\label{hp:ind2}
\tilde{Z}^j_n < Z_{t_i+\tau^\ast_j+n}, \qquad
\tilde{Y}^j_{2,n} \geq Y_{2,t_i+\tau^\ast_j+n},\qquad
\tilde{Y}^j_{1,n} < Y_{1,t_i+\tau^\ast_j+n}.
\end{equation}
In other words, we will show, provided that $t_i>\tau_{\rho_1^{\prime}}$, that for each $j\geq 1$ the process $\tilde{Z}^j_n$ is always
dominated by the original process $Z_{t_i+\tau^\ast_j+n}$, as long as $Z_{t_i+\tau^\ast_j+n}$ is dominated by $\rho_1^{\prime}$
(i.e. for $n\leq {t^\ast_{j+1}}-{\tau^\ast_j}$).
By construction we have that
\[\tilde{Z}^j_0 = \frac{d+u}{2} < u < Z_{t_i+\tau^\ast_j}, \qquad
\tilde{Y}^j_{1,0} = Y_{1,t_i+\tau^\ast_j}\]
which immediately implies $\tilde{Y}^j_{2,0} > Y_{2,t_i+\tau^\ast_j}$.
To this end, we assume~\eqref{hp:ind2} by induction hypothesis.
First, we will show that $\tilde{Y}^j_{2,n+1} > Y_{2,t_i+\tau^\ast_j+n+1}$.
Since from~\eqref{hp:ind2} $\tilde{Z}^j_n < Z_{t_i+\tau^\ast_j+n}$ for $n\leq {t^\ast_{j+1}}-{\tau^\ast_j}$, by construction we obtain that
$$\tilde{X}^j_{n+1} = \ind_{\{\tilde{U}^j_n<\tilde{Z}^j_n\}} \leq
\ind_{\{U_{t_i+\tau^\ast_j+n+1}<Z_{t_i+\tau^\ast_j+n}\}} =X_{t_i+\tau^\ast_j+n+1}.$$
As a consequence, since $W_n\leq1$ for any $n\geq1$, we have that
\[\begin{aligned}
\left(Y_{2,t_i+\tau^\ast_j+n+1}-Y_{2,t_i+\tau^\ast_j+n}\right)\ &&=&\ \left(1-X_{t_i+\tau^\ast_j+n+1}\right)D_{2,t_i+\tau^\ast_j+n+1}
W_{2,t_i+\tau^\ast_j+n}\\
&&\leq&\ (1-\tilde{X}^j_{n+1})\tilde{D}^j_{2,n+1}\\
&&=&\ \left(\tilde{Y}^j_{2,n+1}-\tilde{Y}^j_{2,n}\right),
\end{aligned}\]
which, using hypothesis~\eqref{hp:ind2}, implies $\tilde{Y}^j_{2,n+1} > Y_{2,t_i+\tau^\ast_j+n+1}$.
Similarly, we now show that $\tilde{Y}^j_{1,n+1} \leq Y_{1,t_i+\tau^\ast_j+n+1}$. We have
$$\left(Y_{1,t_i+\tau^\ast_j+n+1}-Y_{1,t_i+\tau^\ast_j+n}\right)\ =\
X_{t_i+\tau^\ast_j+n+1} D_{1,t_i+\tau^\ast_j+n+1}
W_{1,t_i+\tau^\ast_j+n}.$$
From (i) we have that, as long as $Z$ remains below $\rho_1^{\prime}$, $Z$ is also above the process $\hat{\rho}_{1,n}$.
Since we consider the behavior of $Z_{t_i+\tau^\ast_j+n}$ when it is below $\rho_1^{\prime}$, i.e. $n\leq {\tau^\ast_{j+1}}-{t^\ast_j}$,
we have that $W_{1,t_i+\tau^\ast_j+n}=1$. Thus,
$$\left(Y_{1,t_i+\tau^\ast_j+n+1}-Y_{1,t_i+\tau^\ast_j+n}\right)\ \geq\ \tilde{X}^j_{n+1}\tilde{D}^j_{1,n+1}\ =\
\left(\tilde{Y}^j_{1,n+1}-\tilde{Y}^j_{1,n}\right),$$
which using hypothesis~\eqref{hp:ind2} implies $\tilde{Y}^j_{1,n+1} \leq Y_{1,t_i+\tau^\ast_j+n+1}$.
Thus, we have shown that, on the set $\{\tau_{\rho_1^{\prime}}\leq n_{\epsilon}\}$, for any $n\leq {t^\ast_{j+1}}-{\tau^\ast_j}$,
$\tilde{Z}^j_{n+1} < Z_{t_i+\tau^\ast_j+{n+1}}$, $\tilde{Y}^j_{1,n+1} \leq Y_{1,t_i+\tau^\ast_j+{n+1}}$ and
$\tilde{Y}^j_{2,n+1} > Y_{2,t_i+\tau^\ast_j+{n+1}}$ hold.\\

Now, for any $j\geq1$, let $T_j$ be the stopping time for $\tilde{Z}^j_{n}$
to exit from $\left(d,u\right)$, i.e.:
\[\begin{aligned}
T_j & =
\begin{cases}
\inf \{R_3\}\ &\ \text{if } R_3 \neq\emptyset;\\
+ \infty\ &\ \text{otherwise},
\end{cases}
\end{aligned}\]
where $R_3:=\{n\geq1: \tilde{Z}^j_{n}\leq d\text{ or }\tilde{Z}^j_{n}\geq u\}$.
Note that, on the set $\{\tau_{\rho_1^{\prime}}\leq n_\epsilon\}$,
\[\begin{aligned}
\left\{\tau_i<\infty\right\}\ =\ \left\{\inf_{n\geq1}\left\{Z_{t_i+n}\right\}<d\right\}\ &&\subset&\
\left\{\cup_{j: \tau^\ast_{j}\leq n}\left\{\inf_{n\geq1}\left\{\tilde{Z}^j_{n-\tau^\ast_{j}}\right\}<d\right\}\right\}\\
&&\subset&\ \left\{\cup_{j=0}^{\infty}\left\{T_j<\infty\right\}\right\}.
\end{aligned}\]
Hence, by denoting $\bm{P_{i}}\left(\cdot\right)=\bm{P}\left(\cdot|t_i<\infty\right)$ and
$\bm{E_{i}}\left[\cdot\right]=\bm{E}\left[\cdot|t_i<\infty\right]$, we have that
\[\begin{aligned}
\bm{P_{i}}\left(\{\tau_i<\infty\}\cap\{\tau_{\rho_1^{\prime}}\leq n_{\epsilon}\}\right)\ &&\leq&\
\bm{P_{i}}\left(\left\{\cup_{j=0}^{\infty}\left\{T_j<\infty\right\}\right\}\cap\{\tau_{\rho_1^{\prime}}\leq n_{\epsilon}\}\right)\\
&&\leq&\ \sum_{j=0}^{\infty}\bm{P_{i}}\left(\left\{T_j<\infty\right\}\cap\{\tau_{\rho_1^{\prime}}\leq n_{\epsilon}\}\right),
\end{aligned}\]
and, by setting $h=\frac{u-d}{2}$, each term of the series is less or equal than
$$\bm{P_{i}}\left(\left\{\sup_{n\geq1}| \tilde{Z}_{n}^j -
\tilde{Z}_{0}^j| \geq  h\right\}\cap\{\tau_{\rho_1^{\prime}}\leq n_{\epsilon}\}\right)\ \leq\
\bm{P_{i}}\left(\sup_{n\geq1}| \tilde{Z}^j_{n} - \tilde{Z}^j_{0}| \geq  h\right).$$
Note that
$\{\tilde{Z}^j_{n};n\geq1\}$
is the proportion of red balls in an RRU model with same reinforcement means.
Then, by using Lemma~\ref{lem:same_mean} we obtain
\[\begin{aligned}
\bm{P_{i}}\left(\sup_{n\geq1}| \tilde{Z}^j_{n} - \tilde{Z}^j_{0}| \geq  h\right)\ &=&&\
\bm{E_i}\left[ \bm{P}\left(\left.\left\{\sup_{n\geq1}| \tilde{Z}^j_{n} - \tilde{Z}^j_{0}| \geq  h\right\}\right|\mathcal{F}_{\tau_i+t^\ast_{j}}\right) \right]\\
&\leq&&\
\bm{E_i}\left[\frac{b}{{Y}_{t^\ast_{j}}}\right]\left(\frac{4}{h^2}+\frac{2}{h}\right).
\end{aligned}\]
Moreover, by using Lemma~\ref{lem:Y_geometric_increasing}, the right hand side can be expressed as
$$\bm{E_i}\left[\frac{b}{Y_{t_i}}\right]\left(\frac{\rho_1^{\prime}\left(1-\gamma\right)}{\gamma\left(1-\rho_1^{\prime}\right)}\right)^{j}
\left(\frac{4}{h^2}+\frac{2}{h}\right).$$
Since by Lemma~\ref{lem:urn_not_stop} $Y_n$ converges a.s. to infinity, and since $\tau_i\rightarrow\infty$ a.s. because $\tau_i\geq i$,
we have that $\bm{E_i}\left[Y_{t_i}^{-1}\right]$ tends to zero as $i$ increases.
As a consequence, we can choose an integer $i$ large enough such that
$$\bm{E_i}\left[\frac{b}{Y_{t_i}}\right]\left(\frac{4}{h^2}+\frac{2}{h}\right)\left(\frac{1-\rho_1^{\prime}}{1-\rho_1^{\prime}/\gamma}\right)\ <\ \frac{1}{2},$$
which by setting $\phi=1/2+\epsilon$ implies~\eqref{eq:goal_a_s_convergence}, i.e.
$$\bm{P}\left(t_{i+1}<\infty|t_i<\infty\right)\ \leq\ \phi\ <\ 1.$$
This concludes the proof.
\endproof

\proof[Lemma~\ref{lem:harmonic_as_convergence}]
We divide the proof in two parts:
\begin{itemize}
\item[(i)] $m_1\neq m_2$ and $0<\rho_2<\rho_1<1$;
\item[(ii)] $m_1= m_2$ and $0\leq\rho_2<\rho_1\leq1$, on the set $\{Z_{\infty}\neq\{0,1\}\}$;
\end{itemize}

For part (i), assume $m_1>m_2$, since the proof in the case $m_1<m_2$ is completely analogous.
In this case $m^{*}=m_2$ and, by using Theorem~\ref{thm:convergence_as}, we have $Z_n\stackrel{a.s.}{\rightarrow}\rho_1$;
thus, since $\hat{\rho}_{2,n}\stackrel{a.s.}{\rightarrow}\rho_2$ and $\rho_1>\rho_2$,
denoting by $\tau\in\mathbb{N}$ the last time $Z_n$ crosses $\hat{\rho}_{2,n}$, i.e.
$\tau:=\sup\{k\geq1,Z_k<\hat{\rho}_{2,k} \}$, we have that $\bm{P}(\tau<\infty)=1$.
Then, since $\{\tau\leq n\}\subset\{W_{2,k}=1,\forall k\geq n\}$,
we use the following decomposition, on the set $\{\tau\leq n\}$,
$$\frac{Y_{2,n}}{n}\ =\ \frac{1}{n}\sum_{i=1}^n(1-X_i)D_{2,i}W_{2,i-1}\ =\ \mathcal{W}_{0,n}+\mathcal{W}_{1,n},$$
where
\[\begin{aligned}\mathcal{W}_{0,n}\ &&:=&\ \frac{1}{n}\sum_{i=1}^{\tau}(1-X_i)D_{2,i}(W_{2,i-1}-1),\\
\mathcal{W}_{1,n}\ &&:=&\ \frac{1}{n}\sum_{i=\tau}^{n}(1-X_i)D_{2,i}\end{aligned}.\]
Since $\bm{P}(\tau<\infty)=1$, we have $\mathcal{W}_{0,n}\stackrel{a.s.}{\rightarrow}0$, while
since
$$\bm{E}[(1-X_i)D_{2,i}|\mathcal{F}_{i-1}]\ =\ (1-Z_{i-1})m_2\ \stackrel{a.s.}{\rightarrow}\ (1-Z_{\infty})m_2,$$
we have that $\mathcal{W}_{1,n}\stackrel{a.s.}{\rightarrow}(1-Z_{\infty})m_2$.
Finally, since $Y_{n}=(1-Z_n)^{-1}Y_{2,n}$, we have $\frac{Y_{n}}{n}\ \stackrel{a.s.}{\rightarrow}\ m_2=m^{*}$.\\

For part (ii), since $m_1=m_2=m$, by using Theorem~\ref{thm:convergence_as} we have $Z_n\stackrel{a.s.}{\rightarrow}Z_{\infty}\in [\rho_2,\rho_1]$;
then, on the set $\{Z_{\infty}\in(0,1)\}$, we can follow the arguments of part (i), so obtaining
$$\frac{Y_{2,n}}{n}\stackrel{a.s.}{\rightarrow}(1-Z_{\infty})m,\ \qquad\ \frac{Y_{1,n}}{n}\stackrel{a.s.}{\rightarrow}Z_{\infty}m.$$
Thus, $\frac{Y_{n}}{n}=\frac{Y_{1,n}}{n}+\frac{Y_{2,n}}{n}\stackrel{a.s.}{\rightarrow}m.$
\endproof

The proof of Lemma~\ref{lem:no_atoms} is based on comparison arguments between the ARRU and RRU model.
Specifically, for any $n_0\geq1$, we consider an RRU process $\{\tilde{Z}_k(n_0);k\geq0\}$ coupled with the ARRU process $\{Z_{n_0+k};k\geq0\}$ as follows:
the initial composition is $(\tilde{Y}_{1,0}(n_0),\tilde{Y}_{2,0}(n_0))=(Y_{1,n_0},Y_{2,n_0})$ and for any $k\geq1$
\begin{equation}\label{def:dynamics_RRU}\left\{
\begin{array}{l}
\tilde{Y}_{1,k}(n_0)=\tilde{Y}_{1,k-1}(n_0)+\tilde{X}_{1,k}(n_0)D_{1,k}\\
\tilde{Y}_{2,k}(n_0)=\tilde{Y}_{2,k-1}(n_0)+\left(1-\tilde{X}_k(n_0)\right)D_{2,k},
\end{array}
\right.,\end{equation}
where $\tilde{X}_k(n_0)=\ind_{\{U_k\leq\tilde{Z}_{k-1}(n_0)\}}$.
The relation between $\tilde{Z}_k(n_0)$ and $Z_{n_0+k}$ required in the proof of Lemma~\ref{lem:no_atoms} is expressed in the following result.
\begin{lem}\label{lem:relation_RRU_ARRU_CLT_1}
For any $n_0,n_1\geq1$, we have that
\begin{equation}\label{eq:relation_RRU_ARRU_1}
\cap_{k=1}^{n_1}\{ \hat{\rho}_{2,n_0+k}\ \leq\ Z_{n_0+k}\ \leq\ \hat{\rho}_{1,n_0+k}\ \}\ \subset\
\cap_{k=1}^{n_1}\{\ Z_{n_0+k}=\tilde{Z}_k(n_0)\ \}.\end{equation}
\end{lem}
\proof
First, consider the dynamics of the RRU process $\{\tilde{Z}_k(n_0);k\geq0\}$ expressed in~\eqref{def:dynamics_RRU}
and the dynamics of the ARRU process $\{Z_{n_0+k};k\geq0\}$ expressed as follows:
\begin{equation}\label{def:dynamics_ARRU}\left\{
\begin{array}{l}
Y_{1,n_0+k}=Y_{1,n_0+k-1}+X_{1,n_0+k}D_{1,n_0+k}W_{1,n_0+k-1}\\
Y_{2,n_0+k}=Y_{2,n_0+k-1}+\left(1-X_{1,n_0+k}\right)D_{2,n_0+k}W_{2,n_0+k-1},
\end{array}
\right.\end{equation}
where $X_{n_0+k}=\ind_{\{U_k\leq Z_{n_0+k-1}\}}$.
Hence,~\eqref{eq:relation_RRU_ARRU_1} follows by noticing that for any $1\leq k\leq n_1$
$$\left\{ \hat{\rho}_{2,n_0+k}\ \leq\ Z_{n_0+k}\ \leq\ \hat{\rho}_{1,n_0+k}\ \right\}\ \subset\ \{W_{1,n_0+k-1}=W_{2,n_0+k-1}=1\}.$$
\endproof

\proof[Lemma~\ref{lem:no_atoms}]
The proof is structured as follows:
we assume there exist $x\in(\rho_2,\rho_1)$ and $p>0$ such that $\bm{P}(Z_{\infty}=x)=p$
and we show that this assumption leads to a contradiction.
To this end, fix $\epsilon>0$ such that $\rho_2<x-\epsilon<x+\epsilon<\rho_1$ and denote by $\tau\in\mathbb{N}$
the last time $Z_n$ exceeds $I_{\epsilon}:=(x-\epsilon,x+\epsilon)$: formally,
\begin{equation*}
\begin{aligned}
\tau & =
\begin{cases}
\sup\{k>1 :Z_k\notin I_{\epsilon},\}\ \  & \text{if }\{k>1 :Z_k\notin I_{\epsilon}\} \neq\emptyset;
\\
- \infty & \text{otherwise}.
\end{cases}
\end{aligned}
\end{equation*}
Since $\{Z_{\infty}=x\}\subset\{\tau<\infty\}$ and by~\eqref{eq:rho_convergence_as} $\hat{\rho}_{j,n}\stackrel{a.s.}{\rightarrow}\rho_j\notin I_{\epsilon}$, $j\in\{1,2\}$, there exists an integer $k_0\in\mathbb{N}$ such that,
\begin{equation}\label{eq:conradiction_p_2}
\bm{P}\left(\ \{\hat{\rho}_{j,n}\notin I_{\epsilon},\forall n\geq k_0\}\cap\{\tau\leq k_0\}\cap\{Z_{\infty}=x\}\ \right)\ \geq\ \frac{p}{2}.
\end{equation}
Now, by using Lemma~\ref{lem:relation_RRU_ARRU_CLT_1}, we have that
$$\{\hat{\rho}_{j,n}\notin I_{\epsilon},\forall n\geq k_0\}\cap\{\tau\leq k_0\}\ \subset\ \left\{Z_{k_0+n}=\tilde{Z_n}(k_0),\forall n\geq k_0\right\},$$
and hence~\eqref{eq:conradiction_p_2} is equivalent to
$$\bm{P}\left(\ \{\hat{\rho}_{j,n}\notin I_{\epsilon},\forall n\geq k_0\}\cap\{\tau\leq k_0\}\cap\{\tilde{Z}_{\infty}(k_0)=x\}\ \right)\ \geq\ \frac{p}{2}.$$
Finally, the contradiction follows by noticing that by Theorem 2 in~\cite{Aletti.et.al.09}, for RRU model, we have $\bm{P}(\tilde{Z}_{\infty}(k_0)=x)=0$.\\
\endproof

\section{Proofs of limit distribution of the proportion of sampled balls}   \label{section_CLT}

We start by presenting the limit distribution of the proportion of sampled balls for the RRU model.

\proof[Theorem~\ref{thm:CLT_RRU}]
Note that
$$\sqrt{n}\left(\frac{N_{1n}}{n}-Z_{\infty}\right)\ =\ T_{1n}\ +\ T_{2n},$$
where
$$T_{1n}:=n^{-1/2}\left(N_{1n}-\sum_{i=1}^nZ_{i-1}\right),\ \ \ \ T_{2n}:=n^{-1/2}\sum_{i=1}^n\left(Z_{i-1}-Z_{\infty}\right).$$
Now, calling $\Delta Z_j=Z_j-Z_{j-1}$ and $(j\wedge n):=\min\{j,n\}$, we have that
\[\begin{aligned}
T_{2n}\ &&=&\ n^{-1/2}\sum_{i=1}^n\sum_{j=i}^{\infty}(-\Delta Z_j)\ =\ -n^{-1/2}\sum_{j=1}^{\infty}\sum_{i=1}^{j\wedge n}\Delta Z_j\\
&&=&\ -n^{-1/2}\sum_{j=1}^{\infty}(j\wedge n)\Delta Z_j\ =\ -(T_{3n}\ +\ T_{4n}),
\end{aligned}\]
where, since $(j\wedge n)=n$ for all $j\geq n+1$,  we have
$$T_{3n}:=n^{-1/2}\sum_{j=1}^{n}j\Delta Z_j,\ \qquad \ T_{4n}:=n^{1/2}(Z_{\infty}-Z_n).$$
Now, by using the Doob's decomposition $\Delta Z_j=\Delta M_j+\Delta A_j$ (see~\cite{Durrett.11}),
where $\bm{E}[\Delta M_j|\mathcal{F}_{j-1}]=0$ and
$A_j\in\mathcal{F}_{j-1}$, we have $T_{3n}\ =\ T_{5n}\ +\ T_{6n}$, where
$$T_{5n}:=n^{-1/2}\sum_{j=1}^{n}j\Delta M_j,\ \qquad \ T_{6n}:=n^{-1/2}\sum_{j=1}^{n}j\Delta A_j.$$
Then, recalling that
$$\sqrt{n}\left(\frac{N_{1n}}{n}-Z_{\infty}\right)\ =\ T_{1n}\ -\ T_{4n}\ -\ T_{5n}\ -\ T_{6n},$$
the limit distribution is established by proving the following results:
\begin{itemize}
\item[(a)] $T_{4n}|\mathcal{F}_{n}\stackrel{d}{\rightarrow}\mathcal{N}(0,\Sigma_{a})$ (stably), where $\Sigma_{a}=Z_{\infty}(1-Z_{\infty})(1+\frac{\bar{\Sigma}}{m^2})$;
\item[(b)] $T_{6n}\stackrel{p}{\rightarrow}0$;
\item[(c)] $(T_{1n}-T_{5})\stackrel{d}{\rightarrow}\mathcal{N}(0,\Sigma_{c})$ (stably), where
$\Sigma_{c}=Z_{\infty}(1-Z_{\infty})\frac{\bar{\Sigma}}{m^2}$;
\item[(d)] $T_{4n}+(T_{1n}-T_{5n})\stackrel{d}{\rightarrow}\mathcal{N}(0,\Sigma_a+\Sigma_{c})$ (stably).
\end{itemize}

Part (a) follows from Theorem 1 in Aletti \textit{et al.} (see~\cite{Aletti.et.al.09})
and Crimaldi \textit{et al.} (2007) and Crimaldi (2009) (see~\cite{Crimaldi.et.al.07,Crimaldi.09}).

For part (b), by using Lemma~\ref{lem:mul_pag_sec}, for any $j\geq0$, we have that
$$\Delta A_j\ =\ \bm{E}\left[\Delta Z_j|\mathcal{F}_{j-1}\right]\ =\ Z_{j-1}(1-Z_{j-1})B_{j-1},$$
with $W_{1,j-1}=W_{2,j-1}=1$ (since for any $j\geq1$ since the process is an RRU model).
By using Lemma 2 in~\cite{Aletti.et.al.09}, we have $|B_{j-1}|<c_1Y^{-2}_{j-1}$ a.s. for some constant $c_1>0$, and hence
$$T_{6n}\ \leq\ n^{-1/2}\sum_{j=1}^{n}j|\Delta A_j|\ \leq\ c_1n^{-1/2}\sum_{j=1}^{n}jY^{-2}_{j-1};$$
in addition, by using Lemma 3 in~\cite{Aletti.et.al.09}, we have $\bm{E}[Y^{-2}_{j-1}]\leq c_2(j-1)^{-2}$ for some constant $c_2>0$ and hence
$$\bm{E}[T_{6n}]\ \ \leq\ c_1c_2n^{-1/2}\sum_{j=1}^{n}j(j-1)^{-2}\ =\ O\left(n^{-1/2}\log(n)\right).$$
Thus, (b) follows

For part (c), let $T_{1n}-T_{5n}=\sum_{j=1}^n\Delta S_{jn}$ where
$$\Delta S_{jn}\ :=\ n^{-1/2}(X_j-Z_{j-1}-j\Delta M_j).$$
Since $(T_{1n}-T_{5n})$ is a martingale with respect to the filtration $\{\mathcal{F}_n;n\geq1\}$, we apply the Martingale CLT (MCLT)
after establishing the following conditions (see Theorem 3.2 in~\cite{Hall.et.al.80}):
\begin{itemize}
\item[(i)] $\max_{1\leq j\leq n}|\Delta S_{jn}|\stackrel{p}{\rightarrow}0$;
\item[(ii)] $\sup_{n\geq1}\bm{E}[\max_{1\leq j\leq n}(\Delta S_{jn})^2]<\infty$;
\item[(iii)] $\sum_{j=1}^n\bm{E}[(\Delta S_{jn})^2|\mathcal{F}_{j-1}]\stackrel{p}{\rightarrow}\Sigma_c$.
\end{itemize}
For part (i), since $|X_j-Z_{j-1}|\leq1$ a.s. and $\Delta M_j=(\Delta Z_j-\Delta A_j)$, we have that
$$|\Delta S_{jn}|\ \leq\ n^{-1/2}(|X_j-Z_{j-1}|+|j\Delta M_j|)\ \leq\ n^{-1/2}(1+|j(\Delta Z_j-\Delta A_j)|).$$
Now, since $|\Delta Z_j|<bY^{-1}_{j-1}$ and $|\Delta A_j|<c_1Y^{-2}_{j-1}$ a.s. by Lemma 2 in~\cite{Aletti.et.al.09}, we have
$$|\Delta S_{jn}|\ \leq\ n^{-1/2}(1+bjY^{-1}_{j-1}+c_1jY^{-2}_{j-1})\ \ a.s.$$
Since by Lemma~\ref{lem:harmonic_as_convergence} $(jY^{-1}_{j})\stackrel{a.s.}{\rightarrow}m^{-1}$, we have $\sup_{j\geq1}(jY^{-1}_{j})<\infty$ a.s., and thus
$|\Delta S_{jn}|\stackrel{a.s.}{\rightarrow}0$.

For part (ii), using the relation $\bm{E}[S]=\int_0^{\infty}\bm{P}(S>t)dt$ that holds for any non negative r.v. $S$, we obtain
$$\bm{E}\left[\max_{1\leq j\leq n}(\Delta S_{jn})^2\right]\ \leq\ \sum_{j=1}^n\int_{0}^{\infty}\bm{P}((\Delta S_{jn})^2>t)dt.$$
By applying arguments analogous to part (i), we obtain
\[\begin{aligned}n(\Delta S_{jn})^2\ &&\leq&\ 2\left[(X_j-Z_{j-1})^2+(j\Delta M_j)^2\right]\\
&&\leq&\ 2\left[1+2\left[(j\Delta Z_j)^2+(j\Delta A_j)^2\right]\right]\\
&&\leq&\ 2\left[1+2\left[b^2(jY^{-1}_{j-1})^2+c_1^2(jY^{-2}_{j-1})^2\right]\right].
\end{aligned}\]
Thus, by using Markov's inequality we obtain
\[\begin{aligned}\bm{P}((\Delta S_{jn})^2>t)\ &&\leq&\ \bm{P}\left(C\left(\frac{j}{Y_{j-1}}\right)^2>nt\right)\\
&&\leq&\ \max\left\{\ 1\ ;\ \left(\frac{C}{nt}\right)^2\bm{E}\left[\left(\frac{j}{Y_{j-1}}\right)^4\right]\ \right\}.\end{aligned}\]
Now, since by Lemma 3 in~\cite{Aletti.et.al.09} $\sup_{j\geq1}\bm{E}\left[\left(\frac{j}{Y_{j-1}}\right)^4\right]<\infty$,
it follows that there exists a constant $C$ independent of $j$ such that $\int_0^{\infty}\bm{P}((\Delta S_{jn})^2>t)\leq Cn^{-2}$ and hence
$$\sup_{n\geq1}\bm{E}\left[\max_{1\leq j\leq n}(\Delta S_{jn})^2\right]\ \leq\ \sup_{n\geq1}Cn^{-1}\ \leq\ C.$$

For part (iii), since $\Delta M_j=\Delta Z_j-\Delta A_j$, $\Delta A_j\in \mathcal{F}_{j-1}$ and hence
$\bm{E}[\Delta Z_j\Delta A_j|\mathcal{F}_{j-1}]=(\Delta A_j)^2$, we have the following decomposition
$$\bm{E}[(\Delta S_{jn})^2|\mathcal{F}_{j-1}]\ =\ \frac{1}{n}\bm{E}[Q_j^2|\mathcal{F}_{j-1}]\ +\ \frac{2}{n}(j\Delta A_j)^2,$$
where $Q_j:=(X_j-Z_{j-1}-j\Delta Z_j)$. Since $|\Delta A_j|<c_1Y^{-2}_{j-1}$ a.s. and
by Lemma~\ref{lem:harmonic_as_convergence} $(jY^{-1}_{j})\stackrel{a.s.}{\rightarrow}m^{-1}$, we have that $(j\Delta A_j)^2\stackrel{a.s.}{\rightarrow}0$.
Thus, $\frac{2}{n}\sum_{j=1}^n(j\Delta A_j)^2\stackrel{a.s.}{\rightarrow}0$ and hence (iii) is obtained by establishing
\begin{equation}\label{eq:thesis_iii_tilde}
\sum_{j=1}^n\bm{E}[(\Delta S_{jn})^2|\mathcal{F}_{j-1}]\ =\ \frac{1}{n}\sum_{j=1}^n\bm{E}[Q_j^2|\mathcal{F}_{j-1}]\ \stackrel{p}{\rightarrow}\ \Sigma_{c}.
\end{equation}
To this end, we will show that $\bm{E}[Q_j^2|\mathcal{F}_{j-1}]\stackrel{a.s.}{\rightarrow}\Sigma_c$.
First, note that, since $X_j\in\{0,1\}$, we express $\Delta Z_j$ as follows
$$\Delta Z_j\ =\ X_j\left((1-Z_{j-1})\frac{D_{1,j}}{Y_{j-1}}\right)\ +\ (1-X_j)\left(-Z_{j-1}\frac{D_{2,j}}{Y_{j-1}}\right).$$
As a consequence, we consider $Q_j^2=X_jQ_{j,1}^2+(1-X_j)Q_{j,0}^2$, where, denoting by $M_{j-1}:=Y_{j-1}/j$,
\[\begin{aligned}
Q_{j,1}\ &&:=&\ (1-Z_{j-1})\left(1-\frac{D_{1,j}}{M_{j-1}}\right)\ =\ \left(\frac{1-Z_{j-1}}{M_{j-1}}\right)\left(M_{j-1}-D_{1,j}\right),\\
Q_{j,0}\ &&:=&\ Z_{j-1}\left(-1+\frac{D_{2,j}}{M_{j-1}}\right)\ =\ \left(\frac{Z_{j-1}}{M_{j-1}}\right)\left(-M_{j-1}+D_{2,j}\right).
\end{aligned}\]
Then, since $D_{1,j}$, $D_{2,j}$ and $X_{j}$ are independent conditionally to $\mathcal{F}_{j-1}$ and using
\[\begin{aligned}\bm{E}[\left(M_{j-1}-D_{1,j}\right)^2|\mathcal{F}_{j-1}]\ &&=&\ (M_{j-1}-m)^2\ +\ \sigma_1^2,\\
\bm{E}[\left(-M_{j-1}+D_{2,j}\right)^2|\mathcal{F}_{j-1}]\ &&=&\ (M_{j-1}-m)^2\ +\ \sigma_2^2,\end{aligned}\]
we have that
\[\begin{aligned}
\bm{E}[Q_j^2|\mathcal{F}_{j-1}]\ &&=&\ Z_{j-1}\bm{E}[Q_{j,1}^2|\mathcal{F}_{j-1}]\ +\ (1-Z_{j-1})\bm{E}[Q_{j,0}^2|\mathcal{F}_{j-1}]\\
&&=&\ Z_{j-1}\left(\frac{1-Z_{j-1}}{M_{j-1}}\right)^2\left[(M_{j-1}-m)^2+\sigma_1^2\right]\\
&&+&\ (1-Z_{j-1})\left(\frac{Z_{j-1}}{M_{j-1}}\right)^2\left[(M_{j-1}-m)^2+\sigma_2^2\right].
\end{aligned}\]
Finally, since by Lemma~\ref{lem:harmonic_as_convergence} $M_{j-1}\stackrel{a.s.}{\rightarrow}m$ and by Theorem~\ref{thm:convergence_as}
$Z_{j-1}\stackrel{a.s.}{\rightarrow}Z_{\infty}$, it follows that
$$\sum_{j=1}^n\bm{E}[(\Delta \tilde{S}_{jn})^2|\mathcal{F}_{j-1}]\ \stackrel{a.s.}{\rightarrow}\ \Sigma_c\ =\
Z_{\infty}(1-Z_{\infty})\left(\ \frac{\bar{\Sigma}}{m^2}\ \right).$$

For part (d), the result follows by combining part (a), (c), Crimaldi \textit{et al.} (2007) and Crimaldi (2009) (see~\cite{Crimaldi.et.al.07,Crimaldi.09}),
and by noticing that $(T_{1n}-T_{5})\in\mathcal{F}_n$.
\endproof

We now turn to consider the ARRU model.
The limit distribution for the ARRU model can be obtained by Theorem~\ref{thm:CLT_RRU} on the set of trajectories that do not cross
the thresholds $\hat{\rho}_{1,n}$ and $\hat{\rho}_{2,n}$ i.o., and hence $\{Z_{\infty}\in (\rho_2,\rho_1)\}$.
Since this set is not $\mathcal{F}_n$-measurable, we consider a sequence of sets $\{A_n;n\geq1\}$ such that
$\{Z_n\in A_n,ev.\}=\{Z_{\infty}\in (\rho_2,\rho_1)\}$ a.s.
Specifically, we consider the sequence of sets $\{A_n;n\geq1\}$ defined in~\eqref{def:A_n_statement} as follows
\begin{equation}\label{def:A_n}
A_n\ :=\ \left(\ \rho_2+CY_n^{-\alpha}\ ,\ \rho_1-CY_n^{-\alpha}\ \right),
\end{equation}
where $0<C<\infty$ is a positive constant and $0<\alpha< \frac{1}{2}$.
Consider the partition $\Omega=\mathcal{A}_{1}\cup\mathcal{A}_{2}\cup\mathcal{A}_{3}$, where
\begin{equation}\label{def:mathcal_A_n}
\begin{aligned}\mathcal{A}_{1}\ &&:=&\ \{Z_k \in A_k , ev.\},\\
\mathcal{A}_{2}\ &&:=&\ \{Z_k \in A_k , i.o.\} \cap  \{Z_k \not\in A_k , i.o.\},\\
\mathcal{A}_{3}\ &&:=&\ \{Z_k \notin A_k , ev.\}.\end{aligned}\end{equation}
The following lemma establish the relation between $\mathcal{A}_{j}$, $j\in\{1,2,3\}$, and $Z_{\infty}$.
\begin{lem}\label{lem:partition_St_Stau}
Assume $m_1=m_2=m$ and~\eqref{eq:rho_convergence_as} with $\rho_1>\rho_2$. Then,
\begin{itemize}
\item[(a)] $\mathcal{A}_{1}\ =\ \{Z_{\infty}\in(\rho_2,\rho_1)\}$ a.s.;
\item[(b)] $\bm{P}(\mathcal{A}_{2})\ =\ 0$;
\item[(c)] $\mathcal{A}_{3}\ =\ \{Z_{\infty}\in\{\rho_2,\rho_1\}\}$ a.s.
\end{itemize}
\end{lem}
The proof of Lemma~\ref{lem:partition_St_Stau} is based on comparison arguments between the ARRU and an RRU model presented in Lemma~\ref{lem:relation_RRU_ARRU_CLT_1}.
This relation is possible when only one random threshold modify the dynamics of the ARRU.
For this reason, we fix $\epsilon\in(0,(\rho_1-\rho_2)/2)$ and we introduce the following times
\begin{equation}\begin{aligned}\label{def:T_1_2_RRU}
T_1\ &&:=&\ \sup\left\{\ n\geq1\ :\ Z_n>\min\{\hat{\rho}_{1n};\rho_1-\epsilon\}\ \right\},\\
T_2\ &&:=&\ \sup\left\{\ n\geq1\ :\ Z_n<\max\{\hat{\rho}_{2n};\rho_2+\epsilon\}\ \right\}.
\end{aligned}\end{equation}
Let $\mathcal{T}_1:=\{T_1<\infty\}$ and $\mathcal{T}_2:=\{T_2<\infty\}$.
Since $\hat{\rho}_{1n},\hat{\rho}_{2n}$ and $Z_n$ converge a.s., $\bm{P}(\mathcal{T}_1\cup\mathcal{T}_2)=1$.
Then, by comparing the ARRU process with the RRU process defined in~\eqref{def:dynamics_RRU} we have the following result.
\begin{lem}\label{lem:relation_RRU_ARRU_CLT_2}
On the set $\mathcal{T}_1$, for any $n_0,k\geq1$ we have
\begin{equation}\label{eq:relation_RRU_ARRU_2}
\{n_0\geq\mathcal{T}_1\}\ \subset\ \left\{\ \tilde{Z}_k(n_0)\ \leq\ Z_{n_0+k}\ \leq\ \rho_1-\epsilon\ \right\}.
\end{equation}
Analogously, on the set $\mathcal{T}_2$, for any $n_0,k\geq1$ we have
\begin{equation}\label{eq:relation_RRU_ARRU_3}
\{n_0\geq\mathcal{T}_2\}\ \subset\ \left\{\ \rho_2+\epsilon\ \leq\ Z_{n_0+k}\ \leq\ \tilde{Z}_k(n_0)\ \right\}.
\end{equation}
\end{lem}
\proof
Consider the dynamics of the RRU process $\{\tilde{Z}_k(n_0);k\geq0\}$ expressed in~\eqref{def:dynamics_RRU}
and the dynamics of the ARRU process $\{Z_{n_0+k};k\geq0\}$ expressed in~\eqref{def:dynamics_ARRU}.
Then, since $\{n_0\geq\mathcal{T}_1\}\subset\{W_{1,n_0+k-1}=1\}$ and $W_{2,n_0+k-1}\leq1$ we obtain~\eqref{eq:relation_RRU_ARRU_2}.
Analogously, since $\{n_0\geq\mathcal{T}_2\}\subset\{W_{2,n_0+k-1}=1\}$ and $W_{1,n_0+k-1}\leq1$ we have~\eqref{eq:relation_RRU_ARRU_3}.
\endproof

\proof[Lemma~\ref{lem:partition_St_Stau}]
First, let $A:=[\rho_2,\rho_1]$, $t_0=0$ and define for every $j\geq1$
\begin{equation*}
\begin{aligned}
\tau_j & =
\begin{cases}
\inf\{k>t_{j-1} :Z_k\in A_k\}\ \  & \text{if }\{k>t_{j-1} :Z_k\in A_k\} \neq\emptyset;
\\
+ \infty & \text{otherwise}.
\end{cases}
\\
t_j & =
\begin{cases}
\inf\{k>\tau_{j}:\tilde{Z}_{k-\tau_{j}}(\tau_{j})\notin A\} & \text{if }\{k>\tau_{j}:\tilde{Z}_{k-\tau_{j}}(\tau_{j})\notin A\} \neq\emptyset;
\\
+ \infty & \text{otherwise}.
\end{cases}
\end{aligned}
\end{equation*}
Denoting by $T_0$ the last finite time in $\{t_j,\tau_j,j\geq1\}$,
we have the following partition $\Omega=S_t\cup S_{\infty}\cup S_{\tau}$, where
\[\begin{aligned}
S_t\ &&:=&\ \{T_0\in\{t_j,j\geq1\}\}\ =\ \cap_{k\geq T_0}\{Z_k\notin A_k\},\\
S_{\infty}\ &&:=&\ \{T_0=\infty\},\\
S_{\tau}\ &&:=&\ \{T_0\in\{\tau_j,j\geq1\}\}\ =\ \cap_{k\geq T_0}\{\tilde{Z}_{k-T_0}(T_0)\in (\rho_2,\rho_1)\}.
\end{aligned}\]
Thus, we establish the following result:
\begin{itemize}
\item[(i)] $\bm{P}(S_{\infty})=0$,
\item[(ii)] $S_{\tau}\ \subset\ \mathcal{A}_1$, and
\item[(iii)] $S_{\tau}\ \subset\ \{Z_{\infty}\in(\rho_2,\rho_1)\}$.
\end{itemize}

For part (i), this result is obtained by establishing that there exists $i_0\geq 1$ such that, for any $i\geq i_0$,
$$\bm{P}(t_i<\infty|\tau_i<\infty)\ \leq\ \frac{1}{2}.$$
To see this, we recall that by Lemma~\ref{lem:same_mean} we have, for any $h\in(0,1)$,
$$\bm{P}\left( \sup_{k\geq1}|\tilde{Z}_k-\tilde{Z}_0| \geq h \right)\ \leq\ \frac{b}{Y_0}\left(\frac{4}{h^2}+\frac{2}{h}\right)\ \leq\ \frac{6b}{Y_0}h^{-2}.$$
Thus, by using Lemma~\ref{lem:same_mean} with $h=C(\tilde{Y}_0(\tau_{j}))^{-\alpha}$ we obtain
\[\begin{aligned}
\bm{P}(t_i<\infty|\tau_i<\infty)\ &&=&\ \bm{P}\left(\ \cup_{k\geq 1}\tilde{Z}_{k}(\tau_i)\notin [\rho_2,\rho_1]\ \big|\tau_i<\infty\right)\\
&&\leq&\ \bm{P}\left(\ \sup_{k\geq1}|\tilde{Z}_{k}(\tau_{j})-\tilde{Z}_0(\tau_{j})|> C(\tilde{Y}_0(\tau_{j}))^{-\alpha}\ \big|\tau_i<\infty\right)\\
&&\leq&\ \bm{E}\left[\ \left(\frac{6b}{\tilde{Y}_0(\tau_{j})}\right)\left(C(\tilde{Y}_0(\tau_{j}))^{-\alpha}\right)^{-2}\ \big|\tau_i<\infty\right]\\
&&=&\ \frac{6b}{C^2}\bm{E}\left[\ (\tilde{Y}_0(\tau_{j}))^{2\alpha-1}\ \big|\tau_i<\infty\right],
\end{aligned}\]
and hence the result follows by recalling that $0<\alpha<\frac{1}{2}$ and $\tilde{Y}_0(\tau_{j})=Y_{\tau_{j}}\geq Y_0+ja$ a.s.
For part (ii), by Lemma~\ref{lem:relation_RRU_ARRU_CLT_2}, we have that
$$S_{\tau}\cap\mathcal{T}_1\ \subset\ \cap_{k\geq T_0}\{\tilde{Z}_{k-T_0}(T_0)\ \leq\ Z_{k}\ \leq\ \rho_1-\epsilon\},$$
$$S_{\tau}\cap\mathcal{T}_2\ \subset\ \cap_{k\geq T_0}\{\rho_2+\epsilon\ \leq\ Z_{k}\ \leq\ \tilde{Z}_{k-T_0}(T_0)\}.$$
Thus, the result follows by $\bm{P}(\mathcal{T}_1\cup\mathcal{T}_2)=1$ and $\tilde{Z}_{k-T_0}(T_0)\stackrel{a.s.}{\rightarrow}\tilde{Z}_{\infty}(T_0)\in(\rho_2,\rho_1)$.
For part (iii), from part (ii) we have that
$$S_{\tau}\ \subset\ \{\min\{\rho_2+\epsilon,\tilde{Z}_{\infty}(T_0)\}\ \leq\ Z_{\infty}\ \leq\ \max\{\rho_1-\epsilon,\tilde{Z}_{\infty}(T_0)\}\};$$
thus, the result follows by noticing that
$$\left(\ \min\{\rho_2+\epsilon,\tilde{Z}_{\infty}(T_0)\},\max\{\rho_1-\epsilon,\tilde{Z}_{\infty}(T_0)\}\ \right)\ \subset\ (\rho_2,\rho_1).$$
Now, to complete the proof of Lemma~\ref{lem:partition_St_Stau},
we notice that from (i), (ii) and $\{\mathcal{A}_3=S_t\}$, it follows that $\bm{P}(\mathcal{A}_2)=0$ and $\{S_{\tau}=\mathcal{A}_1\}$.
Then, combining (iii) and $\mathcal{A}_3\subset \{Z_{\infty}\in\{\rho_2,\rho_1\}\}$, we obtain the result.
\endproof

We now present the proof of the limit distribution of the proportion of sampled balls for the ARRU model.

\proof[Theorem~\ref{thm:CLT_ARRU}]
First, take the sets $\mathcal{A}_1$, $\mathcal{A}_2$ and $\mathcal{A}_3$ defined in~\eqref{def:mathcal_A_n}.
Note that, since $\mathcal{A}_1=\underline{\lim}_n \{Z_n\in A_n\}$ and $\mathcal{A}^c_3=\overline{\lim}_n \{Z_n\in A_n\}$,
by Lemma~\ref{lem:partition_St_Stau} we have
$$\underline{\lim}_n \{Z_n\in A_n\}\ =\ \overline{\lim}_n \{Z_n\in A_n\}\ =\ \{Z_{\infty}\in (\rho_2,\rho_1)\}.$$
Then, the proof is based on applying Theorem~\ref{thm:CLT_RRU} to the ARRU model.
To this end, consider the decomposition
$\{Z_n\in A_n\}=\mathcal{A}_{1n}\cup\mathcal{A}_{2n}\cup\mathcal{A}_{3n}$, where $\mathcal{A}_{jn}=\{Z_n\in A_n\}\cap \mathcal{A}_{j}$ for any $j\in\{1,2,3\}$.
Since by using Lemma~\ref{def:A_n} $\bm{P}(\mathcal{A}_{2})=0$, we have $\bm{P}(\mathcal{A}_{2n})=0$ for any $n\geq1$.
Moreover, by definition we have that $\bm{P}(\mathcal{A}_{3n})\rightarrow0$ and $\bm{P}(\mathcal{A}_{1n})\rightarrow \bm{P}(\mathcal{A}_{1})$.
Thus, calling $\mathcal{N}_n:=\sqrt{n}(\frac{N_{1n}}{n}-Z_{\infty})$, we have
$$\lim_{n\rightarrow\infty}\bm{P}\left(\ \mathcal{N}_n\leq x\ ,\ \{Z_n\in A_n\}\ \right)\ =\
\lim_{n\rightarrow\infty}\bm{P}\left(\ \mathcal{N}_n\leq x\ ,\ \mathcal{A}_1\ \right),$$
and since by Lemma~\ref{def:A_n} $\mathcal{A}_{1}=\{Z_{\infty}\in(\rho_2,\rho_1)\}$, this is equivalent to
$$\lim_{n\rightarrow\infty}\bm{P}\left(\ \mathcal{N}_n\leq x\ ,\ \{Z_{\infty}\in(\rho_2,\rho_1)\}\ \right).$$

Now, consider the RRU model $\{\tilde{Z}_{k}(n_0),k\geq1\}$ described in~\eqref{def:dynamics_RRU} coupled with the ARRU model $\{Z_{n_0+k},k\geq1\}$.
By using Lemma~\ref{lem:relation_RRU_ARRU_CLT_1}, for any $n_0\geq1$, we have
$$\cap_{k=n_0}^{\infty}\{ \hat{\rho}_{2,k}\ \leq\ Z_{k}\ \leq\ \hat{\rho}_{1,k}\ \}\ \subset\
\cap_{k=1}^{\infty}\{ Z_{n_0+k}=\tilde{Z}_k(n_0)\ \}.$$
Hence, on this set the ARRU process $Z_{n_0+k}$ is equivalent to the RRU process $\tilde{Z}_k(n_0)$;
thus, we can obtain the limit distribution for the ARRU by applying the limit distribution for the RRU expressed in Theorem~\ref{thm:CLT_RRU} on the set
where the trajectories of the two processes are equivalent.
To this end, define
$$T^{*}\ :=\ \sup\left\{\ k\geq1\ :\ \{Z_k<\hat{\rho}_{2,k}\}\cup\{Z_k>\hat{\rho}_{1,k}\}\ \right\},$$
and note that, for any $n_0\geq 1$,
$$\{T^{*}\leq n_0\}\ \subset\ \cap_{k=1}^{\infty}\{ Z_{n_0+k}=\tilde{Z}_k(n_0)\ \}.$$
Let $\mathcal{S}$ be a r.v. with characteristic function $\bm{E}[\exp(\frac{1}{2}\Sigma t^2)]$.
Thus, by applying Theorem~\ref{thm:CLT_RRU} we have that, for any $n_0\geq 1$ and any set $\mathcal{T}\in\mathcal{F}$,
$$\lim_{n\rightarrow\infty}\bm{P}\left(\ \mathcal{N}_n\leq x\ ,\ \mathcal{T}\cap\{T^{*}\leq n_0\}\ \right)\ =\
\bm{P}(\ \mathcal{S}\leq x\ ,\ \mathcal{T}\cap \{T^{*}\leq n_0\}\ ).$$
Now, since $\{Z_{\infty}\in(\rho_2,\rho_1)\}\subset\{T^{*}<\infty\}$,
we have
$$\lim_{n_0\rightarrow\infty}\bm{P}(\{T^{*}\leq n_0\}\cap\{Z_{\infty}\in(\rho_2,\rho_1)\})\ =\ \bm{P}(Z_{\infty}\in(\rho_2,\rho_1)),$$
which implies that
$$\lim_{n\rightarrow\infty}\bm{P}\left(\ \mathcal{N}_n\leq x,\{Z_{\infty}\in(\rho_2,\rho_1)\}\ \right)\ =\
\bm{P}(\ \mathcal{S}\leq x ,  \{Z_{\infty}\in (\rho_2,\rho_1)\}\ ).$$
This concludes the proof.
\endproof

E-mail: {\tt giacomo.aletti@unimi.it}, {\tt andrea.ghiglietti@unimi.it}, {\tt avidyash@gmu.edu}.
\end{document}